\documentclass[12pt]{article}

\begin{document}
\setlength{\unitlength}{20pt}
\def\DottedCircle{
\qbezier[4](0.966,-0.259)(1.04,0)(0.966,0.259)
\qbezier[4](0.966,0.259)(0.897,0.518)(0.707,0.707)
\qbezier[4](0.707,0.707)(0.518,0.897)(0.259,0.966)
\qbezier[4](0.259,0.966)(0,1.04)(-0.259,0.966)
\qbezier[4](-0.259,0.966)(-0.518,0.897)(-0.707,0.707)
\qbezier[4](-0.707,0.707)(-0.897,0.518)(-0.966,0.259)
\qbezier[4](-0.966,0.259)(-1.04,0)(-0.966,-0.259)
\qbezier[4](-0.966,-0.259)(-0.897,-0.518)(-0.707,-0.707)
\qbezier[4](-0.707,-0.707)(-0.518,-0.897)(-0.259,-0.966)
\qbezier[4](-0.259,-0.966)(0,-1.04)(0.259,-0.966)
\qbezier[4](0.259,-0.966)(0.518,-0.897)(0.707,-0.707)
\qbezier[4](0.707,-0.707)(0.897,-0.518)(0.966,-0.259)
}
\def\FullCircle{
\thicklines
\put(0,0){\circle{2}}
}
%
%
\def\Endpoint[#1]{
\ifcase#1
\put(1,0){\circle*{0.15}}
\or\put(0.866,0.5){\circle*{0.15}}
\or\put(0.5,0.866){\circle*{0.15}}
\or\put(0,1){\circle*{0.15}}
\or\put(-0.5,0.866){\circle*{0.15}}
\or\put(-0.866,0.5){\circle*{0.15}}
\or\put(-1,0){\circle*{0.15}}
\or\put(-0.866,-0.5){\circle*{0.15}}
\or\put(-0.5,-0.866){\circle*{0.15}}
\or\put(0,-1){\circle*{0.15}}
\or\put(0.5,-0.866){\circle*{0.15}}
\or\put(0.866,-0.5){\circle*{0.15}}
\fi}
%
%
\def\Arc[#1]{
\ifcase#1
\qbezier[25](0.966,-0.259)(1.04,0)(0.966,0.259)
\or
\qbezier[25](0.966,0.259)(0.897,0.518)(0.707,0.707)
\or
\qbezier[25](0.707,0.707)(0.518,0.897)(0.259,0.966)
\or
\qbezier[25](0.259,0.966)(0,1.04)(-0.259,0.966)
\or
\qbezier[25](-0.259,0.966)(-0.518,0.897)(-0.707,0.707)
\or
\qbezier[25](-0.707,0.707)(-0.897,0.518)(-0.966,0.259)
\or
\qbezier[25](-0.966,0.259)(-1.04,0)(-0.966,-0.259)
\or
\qbezier[25](-0.966,-0.259)(-0.897,-0.518)(-0.707,-0.707)
\or
\qbezier[25](-0.707,-0.707)(-0.518,-0.897)(-0.259,-0.966)
\or
\qbezier[25](-0.259,-0.966)(0,-1.04)(0.259,-0.966)
\or
\qbezier[25](0.259,-0.966)(0.518,-0.897)(0.707,-0.707)
\or
\qbezier[25](0.707,-0.707)(0.897,-0.518)(0.966,-0.259)
\fi}
%
%
\def\Chord[#1,#2]{
\thinlines
\ifnum#1>#2\Chord[#2,#1]
\else\ifnum#1<#2
\ifcase#1
\ifcase#2
\or\qbezier(1,0)(0.516,0.138)(0.866,0.5)
\or\qbezier(1,0)(0.45,0.26)(0.5,0.866)
\or\qbezier(1,0)(0.327,0.327)(0,1)
\or\qbezier(1,0)(0.179,0.311)(-0.5,0.866)
\or\qbezier(1,0)(0.0536,0.2)(-0.866,0.5)
\or\put(1, 0){\line(-2, 0){2}}
\or\qbezier(1,0)(0.0536,-0.2)(-0.866,-0.5)
\or\qbezier(1,0)(0.179,-0.311)(-0.5,-0.866)
\or\qbezier(1,0)(0.327,-0.327)(0,-1)
\or\qbezier(1,0)(0.45,-0.26)(0.5,-0.866)
\or\qbezier(1,0)(0.516,-0.138)(0.866,-0.5)
\fi
\or\ifcase#2\or
\or\qbezier(0.866,0.5)(0.378,0.378)(0.5,0.866)
\or\qbezier(0.866,0.5)(0.26,0.45)(0,1)
\or\qbezier(0.866,0.5)(0.12,0.446)(-0.5,0.866)
\or\qbezier(0.866,0.5)(0,0.359)(-0.866,0.5)
\or\qbezier(0.866,0.5)(-0.0536,0.2)(-1,0)
\or\put(0.866, 0.5){\line(-5, -3){1.73}}
\or\qbezier(0.866,0.5)(0.146,-0.146)(-0.5,-0.866)
\or\qbezier(0.866,0.5)(0.311,-0.179)(0,-1)
\or\qbezier(0.866,0.5)(0.446,-0.12)(0.5,-0.866)
\or\qbezier(0.866,0.5)(0.52,0)(0.866,-0.5)
\fi
\or\ifcase#2\or\or
\or\qbezier(0.5,0.866)(0.138,0.516)(0,1)
\or\qbezier(0.5,0.866)(0,0.52)(-0.5,0.866)
\or\qbezier(0.5,0.866)(-0.12,0.446)(-0.866,0.5)
\or\qbezier(0.5,0.866)(-0.179,0.311)(-1,0)
\or\qbezier(0.5,0.866)(-0.146,0.146)(-0.866,-0.5)
\or\put(0.5, 0.866){\line(-3, -5){1}}
\or\qbezier(0.5,0.866)(0.2,-0.0536)(0,-1)
\or\qbezier(0.5,0.866)(0.359,0)(0.5,-0.866)
\or\qbezier(0.5,0.866)(0.446,0.12)(0.866,-0.5)
\fi
\or\ifcase#2\or\or\or
\or\qbezier(0,1.)(-0.138,0.516)(-0.5,0.866)
\or\qbezier(0,1.)(-0.26,0.45)(-0.866,0.5)
\or\qbezier(0,1.)(-0.327,0.327)(-1,0)
\or\qbezier(0,1.)(-0.311,0.179)(-0.866,-0.5)
\or\qbezier(0,1.)(-0.2,0.0536)(-0.5,-0.866)
\or\put(0, 1){\line(0, -2){2}}
\or\qbezier(0,1.)(0.2,0.0536)(0.5,-0.866)
\or\qbezier(0,1.)(0.311,0.179)(0.866,-0.5)
\fi
\or\ifcase#2\or\or\or\or
\or\qbezier(-0.5,0.866)(-0.378,0.378)(-0.866,0.5)
\or\qbezier(-0.5,0.866)(-0.45,0.26)(-1,0)
\or\qbezier(-0.5,0.866)(-0.446,0.12)(-0.866,-0.5)
\or\qbezier(-0.5,0.866)(-0.359,0)(-0.5,-0.866)
\or\qbezier(-0.5,0.866)(-0.2,-0.0536)(0,-1)
\or\put(-0.5, 0.866){\line(3, -5){1}}
\or\qbezier(-0.5,0.866)(0.146,0.146)(0.866,-0.5)
\fi
\or\ifcase#2\or\or\or\or\or
\or\qbezier(-0.866,0.5)(-0.516,0.138)(-1,0)
\or\qbezier(-0.866,0.5)(-0.52,0)(-0.866,-0.5)
\or\qbezier(-0.866,0.5)(-0.446,-0.12)(-0.5,-0.866)
\or\qbezier(-0.866,0.5)(-0.311,-0.179)(0,-1)
\or\qbezier(-0.866,0.5)(-0.146,-0.146)(0.5,-0.866)
\or\put(-0.866, 0.5){\line(5, -3){1.73}}
\fi
\or\ifcase#2\or\or\or\or\or\or
\or\qbezier(-1,0)(-0.516,-0.138)(-0.866,-0.5)
\or\qbezier(-1,0)(-0.45,-0.26)(-0.5,-0.866)
\or\qbezier(-1,0)(-0.327,-0.327)(0,-1)
\or\qbezier(-1,0)(-0.179,-0.311)(0.5,-0.866)
\or\qbezier(-1,0)(-0.0536,-0.2)(0.866,-0.5)
\fi
\or\ifcase#2\or\or\or\or\or\or\or
\or\qbezier(-0.866,-0.5)(-0.378,-0.378)(-0.5,-0.866)
\or\qbezier(-0.866,-0.5)(-0.26,-0.45)(0,-1)
\or\qbezier(-0.866,-0.5)(-0.12,-0.446)(0.5,-0.866)
\or\qbezier(-0.866,-0.5)(0,-0.359)(0.866,-0.5)
\fi
\or\ifcase#2\or\or\or\or\or\or\or\or
\or\qbezier(-0.5,-0.866)(-0.138,-0.516)(0,-1)
\or\qbezier(-0.5,-0.866)(0,-0.52)(0.5,-0.866)
\or\qbezier(-0.5,-0.866)(0.12,-0.446)(0.866,-0.5)
\fi
\or\ifcase#2\or\or\or\or\or\or\or\or\or
\or\qbezier(0,-1.)(0.138,-0.516)(0.5,-0.866)
\or\qbezier(0,-1.)(0.26,-0.45)(0.866,-0.5)
\fi
\or\ifcase#2\or\or\or\or\or\or\or\or\or\or
\or\qbezier(0.5,-0.866)(0.378,-0.378)(0.866,-0.5)
\fi\fi\fi\fi}
%
%
\def\FullChord[#1,#2]{
\Endpoint[#1]
\Endpoint[#2]
\Chord[#1,#2]
}
%
%
\def\EndChord[#1,#2]{
\Endpoint[#1]
\Endpoint[#2]
\Chord[#1,#2]
}
%
%
%
%
%
%
\def\Picture#1{
\begin{picture}(2,1)(-1,-0.167)
#1
\end{picture}
}
%
%
\def\DottedChordDiagram[#1,#2]{
\Picture{\DottedCircle \FullChord[#1,#2]}
}
\let\DCD=\DottedChordDiagram
%
%
%
%

\def\SelfInt{\put(-1,-1){\vector(1,1){2}}\put(-1,1){\vector(1,-1){2}}
\put(0,0){\circle*{0.15}}}
\def\Over{\put(-1,-1){\vector(1,1){2}}\put(-1,1){\vector(1,-1){0.75}}\put(0.25,-0.25){\vector(1,-1){0.75}}}
\def\Under{\put(-1,-1){\vector(1,1){0.75}}\put(0.25,0.25){\vector(1,1){0.75}}\put(-1,1){\vector(1,-1){2}}}

\def\Sterm{
\DottedCircle
\Arc[8]
\Arc[9]
\Arc[10]
\Endpoint[9]
\put(0,-0.5){\circle*{0.15}}
\put(0,-1){\line(0,+1){0.5}}
\put(0,-0.5){\vector(1,1){0.5}}
\put(0,-0.5){\vector(-1,1){0.5}}}
\def\Tterm{
\DottedCircle
\Arc[8]
\Arc[9]
\Arc[10]
\qbezier(-0.5,0)(-0.1,-0.5)(-0.1,-0.97)
\qbezier(+0.5,0)(+0.1,-0.5)(+0.1,-0.97)
\put(-0.1,-0.97){\circle*{0.15}}
\put(0.1,-0.97){\circle*{0.15}}
\put(-0.5,0){\vector(-2,3){0.1}}
\put(0.5,0){\vector(2,3){0.1}}
}
\def\Uterm{
\DottedCircle
\Arc[8]
\Arc[9]
\Arc[10]
\put(-0.1,-0.97){\circle*{0.15}}
\put(0.1,-0.97){\circle*{0.15}}
\qbezier(-0.5,0)(0.1,-0.5)(0.1,-0.97)
\qbezier(+0.5,0)(-0.1,-0.5)(-.1,-0.97)
\put(-0.5,0){\vector(-3,4){0.1}}
\put(0.5,0){\vector(3,4){0.1}}
}

\def\InOut[#1,#2]{
\thinlines
\ifcase#2
\ifcase#1
\qbezier(1,0)(1,0)(0.667,0)
\or\qbezier(1,0)(1,0)(0.206,0.634)
\or\qbezier(1,0)(1,0)(-0.539,0.392)
\or\qbezier(1,0)(1,0)(-0.539,-0.391)
\or\qbezier(1,0)(1,0)(0.206,-0.634)
\or\qbezier(1,0)(1,0)(0.333,0)
\or\qbezier(1,0)(1,0)(-0.167,0.289)
\or\qbezier(1,0)(1,0)(-0.167,-0.289)
\or\qbezier(1,0)(1,0)(0,0)
\fi
\or\ifcase#1
\qbezier(0.866,0.5)(0.866,0.5)(0.667,0)
\or\qbezier(0.866,0.5)(0.866,0.5)(0.206,0.634)
\or\qbezier(0.866,0.5)(0.866,0.5)(-0.539,0.392)
\or\qbezier(0.866,0.5)(0.866,0.5)(-0.539,-0.391)
\or\qbezier(0.866,0.5)(0.866,0.5)(0.206,-0.634)
\or\qbezier(0.866,0.5)(0.866,0.5)(0.333,0)
\or\qbezier(0.866,0.5)(0.866,0.5)(-0.167,0.289)
\or\qbezier(0.866,0.5)(0.866,0.5)(-0.167,-0.289)
\or\qbezier(0.866,0.5)(0.866,0.5)(0,0)
\fi
\or\ifcase#1
\qbezier(0.5,0.866)(0.5,0.866)(0.667,0)
\or\qbezier(0.5,0.866)(0.5,0.866)(0.206,0.634)
\or\qbezier(0.5,0.866)(0.5,0.866)(-0.539,0.392)
\or\qbezier(0.5,0.866)(0.5,0.866)(-0.539,-0.391)
\or\qbezier(0.5,0.866)(0.5,0.866)(0.206,-0.634)
\or\qbezier(0.5,0.866)(0.5,0.866)(0.333,0)
\or\qbezier(0.5,0.866)(0.5,0.866)(-0.167,0.289)
\or\qbezier(0.5,0.866)(0.5,0.866)(-0.167,-0.289)
\or\qbezier(0.5,0.866)(0.5,0.866)(0,0)
\fi
\or\ifcase#1
\qbezier(0,1)(0,1)(0.667,0)
\or\qbezier(0,1)(0,1)(0.206,0.634)
\or\qbezier(0,1)(0,1)(-0.539,0.392)
\or\qbezier(0,1)(0,1)(-0.539,-0.391)
\or\qbezier(0,1)(0,1)(0.206,-0.634)
\or\qbezier(0,1)(0,1)(0.333,0)
\or\qbezier(0,1)(0,1)(-0.167,0.289)
\or\qbezier(0,1)(0,1)(-0.167,-0.289)
\or\qbezier(0,1)(0,1)(0,0)
\fi
\or\ifcase#1
\qbezier(-0.5,0.866)(-0.5,0.866)(0.667,0)
\or\qbezier(-0.5,0.866)(-0.5,0.866)(0.206,0.634)
\or\qbezier(-0.5,0.866)(-0.5,0.866)(-0.539,0.392)
\or\qbezier(-0.5,0.866)(-0.5,0.866)(-0.539,-0.391)
\or\qbezier(-0.5,0.866)(-0.5,0.866)(0.206,-0.634)
\or\qbezier(-0.5,0.866)(-0.5,0.866)(0.333,0)
\or\qbezier(-0.5,0.866)(-0.5,0.866)(-0.167,0.289)
\or\qbezier(-0.5,0.866)(-0.5,0.866)(-0.167,-0.289)
\or\qbezier(-0.5,0.866)(-0.5,0.866)(0,0)
\fi
\or\ifcase#1
\qbezier(-0.866,0.5)(-0.866,0.5)(0.667,0)
\or\qbezier(-0.866,0.5)(-0.866,0.5)(0.206,0.634)
\or\qbezier(-0.866,0.5)(-0.866,0.5)(-0.539,0.392)
\or\qbezier(-0.866,0.5)(-0.866,0.5)(-0.539,-0.391)
\or\qbezier(-0.866,0.5)(-0.866,0.5)(0.206,-0.634)
\or\qbezier(-0.866,0.5)(-0.866,0.5)(0.333,0)
\or\qbezier(-0.866,0.5)(-0.866,0.5)(-0.167,0.289)
\or\qbezier(-0.866,0.5)(-0.866,0.5)(-0.167,-0.289)
\or\qbezier(-0.866,0.5)(-0.866,0.5)(0,0)
\fi
\or\ifcase#1
\qbezier(-1,0)(-1,0)(0.667,0)
\or\qbezier(-1,0)(-1,0)(0.206,0.634)
\or\qbezier(-1,0)(-1,0)(-0.539,0.392)
\or\qbezier(-1,0)(-1,0)(-0.539,-0.391)
\or\qbezier(-1,0)(-1,0)(0.206,-0.634)
\or\qbezier(-1,0)(-1,0)(0.333,0)
\or\qbezier(-1,0)(-1,0)(-0.167,0.289)
\or\qbezier(-1,0)(-1,0)(-0.167,-0.289)
\or\qbezier(-1,0)(-1,0)(0,0)
\fi
\or\ifcase#1
\qbezier(-0.866,-0.5)(-0.866,-0.5)(0.667,0)
\or\qbezier(-0.866,-0.5)(-0.866,-0.5)(0.206,0.634)
\or\qbezier(-0.866,-0.5)(-0.866,-0.5)(-0.539,0.392)
\or\qbezier(-0.866,-0.5)(-0.866,-0.5)(-0.539,-0.391)
\or\qbezier(-0.866,-0.5)(-0.866,-0.5)(0.206,-0.634)
\or\qbezier(-0.866,-0.5)(-0.866,-0.5)(0.333,0)
\or\qbezier(-0.866,-0.5)(-0.866,-0.5)(-0.167,0.289)
\or\qbezier(-0.866,-0.5)(-0.866,-0.5)(-0.167,-0.289)
\or\qbezier(-0.866,-0.5)(-0.866,-0.5)(0,0)
\fi
\or\ifcase#1
\qbezier(-0.5,-0.866)(-0.5,-0.866)(0.667,0)
\or\qbezier(-0.5,-0.866)(-0.5,-0.866)(0.206,0.634)
\or\qbezier(-0.5,-0.866)(-0.5,-0.866)(-0.539,0.392)
\or\qbezier(-0.5,-0.866)(-0.5,-0.866)(-0.539,-0.391)
\or\qbezier(-0.5,-0.866)(-0.5,-0.866)(0.206,-0.634)
\or\qbezier(-0.5,-0.866)(-0.5,-0.866)(0.333,0)
\or\qbezier(-0.5,-0.866)(-0.5,-0.866)(-0.167,0.289)
\or\qbezier(-0.5,-0.866)(-0.5,-0.866)(-0.167,-0.289)
\or\qbezier(-0.5,-0.866)(-0.5,-0.866)(0,0)
\fi
\or\ifcase#1
\qbezier(0,-1)(0,-1)(0.667,0)
\or\qbezier(0,-1)(0,-1)(0.206,0.634)
\or\qbezier(0,-1)(0,-1)(-0.539,0.392)
\or\qbezier(0,-1)(0,-1)(-0.539,-0.391)
\or\qbezier(0,-1)(0,-1)(0.206,-0.634)
\or\qbezier(0,-1)(0,-1)(0.333,0)
\or\qbezier(0,-1)(0,-1)(-0.167,0.289)
\or\qbezier(0,-1)(0,-1)(-0.167,-0.289)
\or\qbezier(0,-1)(0,-1)(0,0)
\fi
\or\ifcase#1
\qbezier(0.5,-0.866)(0.5,-0.866)(0.667,0)
\or\qbezier(0.5,-0.866)(0.5,-0.866)(0.206,0.634)
\or\qbezier(0.5,-0.866)(0.5,-0.866)(-0.539,0.392)
\or\qbezier(0.5,-0.866)(0.5,-0.866)(-0.539,-0.391)
\or\qbezier(0.5,-0.866)(0.5,-0.866)(0.206,-0.634)
\or\qbezier(0.5,-0.866)(0.5,-0.866)(0.333,0)
\or\qbezier(0.5,-0.866)(0.5,-0.866)(-0.167,0.289)
\or\qbezier(0.5,-0.866)(0.5,-0.866)(-0.167,-0.289)
\or\qbezier(0.5,-0.866)(0.5,-0.866)(0,0)
\fi
\or\ifcase#1
\qbezier(0.866,-0.5)(0.866,-0.5)(0.667,0)
\or\qbezier(0.866,-0.5)(0.866,-0.5)(0.206,0.634)
\or\qbezier(0.866,-0.5)(0.866,-0.5)(-0.539,0.392)
\or\qbezier(0.866,-0.5)(0.866,-0.5)(-0.539,-0.391)
\or\qbezier(0.866,-0.5)(0.866,-0.5)(0.206,-0.634)
\or\qbezier(0.866,-0.5)(0.866,-0.5)(0.333,0)
\or\qbezier(0.866,-0.5)(0.866,-0.5)(-0.167,0.289)
\or\qbezier(0.866,-0.5)(0.866,-0.5)(-0.167,-0.289)
\or\qbezier(0.866,-0.5)(0.866,-0.5)(0,0)
\fi
\fi
}

\def\InIn[#1,#2]{
\thinlines
\ifcase#1
\ifcase#2
\qbezier(0.667,0)(0.667,0)(0.667,0)
\or\qbezier(0.667,0)(0.667,0)(0.206,0.634)
\or\qbezier(0.667,0)(0.667,0)(-0.539,0.392)
\or\qbezier(0.667,0)(0.667,0)(-0.539,-0.391)
\or\qbezier(0.667,0)(0.667,0)(0.206,-0.634)
\or\qbezier(0.667,0)(0.667,0)(0.333,0)
\or\qbezier(0.667,0)(0.667,0)(-0.167,0.289)
\or\qbezier(0.667,0)(0.667,0)(-0.167,-0.289)
\or\qbezier(0.667,0)(0.667,0)(0,0)
\fi
\or\ifcase#2
\qbezier(0.206,0.634)(0.206,0.634)(0.667,0)
\or\qbezier(0.206,0.634)(0.206,0.634)(0.206,0.634)
\or\qbezier(0.206,0.634)(0.206,0.634)(-0.539,0.392)
\or\qbezier(0.206,0.634)(0.206,0.634)(-0.539,-0.391)
\or\qbezier(0.206,0.634)(0.206,0.634)(0.206,-0.634)
\or\qbezier(0.206,0.634)(0.206,0.634)(0.333,0)
\or\qbezier(0.206,0.634)(0.206,0.634)(-0.167,0.289)
\or\qbezier(0.206,0.634)(0.206,0.634)(-0.167,-0.289)
\or\qbezier(0.206,0.634)(0.206,0.634)(0,0)
\fi
\or\ifcase#2
\qbezier(-0.539,0.392)(-0.539,0.392)(0.667,0)
\or\qbezier(-0.539,0.392)(-0.539,0.392)(0.206,0.634)
\or\qbezier(-0.539,0.392)(-0.539,0.392)(-0.539,0.392)
\or\qbezier(-0.539,0.392)(-0.539,0.392)(-0.539,-0.391)
\or\qbezier(-0.539,0.392)(-0.539,0.392)(0.206,-0.634)
\or\qbezier(-0.539,0.392)(-0.539,0.392)(0.333,0)
\or\qbezier(-0.539,0.392)(-0.539,0.392)(-0.167,0.289)
\or\qbezier(-0.539,0.392)(-0.539,0.392)(-0.167,-0.289)
\or\qbezier(-0.539,0.392)(-0.539,0.392)(0,0)
\fi
\or\ifcase#2
\qbezier(-0.539,-0.391)(-0.539,-0.391)(0.667,0)
\or\qbezier(-0.539,-0.391)(-0.539,-0.391)(0.206,0.634)
\or\qbezier(-0.539,-0.391)(-0.539,-0.391)(-0.539,0.392)
\or\qbezier(-0.539,-0.391)(-0.539,-0.391)(-0.539,-0.391)
\or\qbezier(-0.539,-0.391)(-0.539,-0.391)(0.206,-0.634)
\or\qbezier(-0.539,-0.391)(-0.539,-0.391)(0.333,0)
\or\qbezier(-0.539,-0.391)(-0.539,-0.391)(-0.167,0.289)
\or\qbezier(-0.539,-0.391)(-0.539,-0.391)(-0.167,-0.289)
\or\qbezier(-0.539,-0.391)(-0.539,-0.391)(0,0)
\fi
\or\ifcase#2
\qbezier(0.206,-0.634)(0.206,-0.634)(0.667,0)
\or\qbezier(0.206,-0.634)(0.206,-0.634)(0.206,0.634)
\or\qbezier(0.206,-0.634)(0.206,-0.634)(-0.539,0.392)
\or\qbezier(0.206,-0.634)(0.206,-0.634)(-0.539,-0.391)
\or\qbezier(0.206,-0.634)(0.206,-0.634)(0.206,-0.634)
\or\qbezier(0.206,-0.634)(0.206,-0.634)(0.333,0)
\or\qbezier(0.206,-0.634)(0.206,-0.634)(-0.167,0.289)
\or\qbezier(0.206,-0.634)(0.206,-0.634)(-0.167,-0.289)
\or\qbezier(0.206,-0.634)(0.206,-0.634)(0,0)
\fi
\or\ifcase#2
\qbezier(0.333,0)(0.333,0)(0.667,0)
\or\qbezier(0.333,0)(0.333,0)(0.206,0.634)
\or\qbezier(0.333,0)(0.333,0)(-0.539,0.392)
\or\qbezier(0.333,0)(0.333,0)(-0.539,-0.391)
\or\qbezier(0.333,0)(0.333,0)(0.206,-0.634)
\or\qbezier(0.333,0)(0.333,0)(0.333,0)
\or\qbezier(0.333,0)(0.333,0)(-0.167,0.289)
\or\qbezier(0.333,0)(0.333,0)(-0.167,-0.289)
\or\qbezier(0.333,0)(0.333,0)(0,0)
\fi
\or\ifcase#2
\qbezier(-0.167,0.289)(-0.167,0.289)(0.667,0)
\or\qbezier(-0.167,0.289)(-0.167,0.289)(0.206,0.634)
\or\qbezier(-0.539,0.392)(-0.167,0.289)(-0.167,0.289)
\or\qbezier(-0.539,-0.391)(-0.167,0.289)(-0.167,0.289)
\or\qbezier(0.206,-0.634)(-0.167,0.289)(-0.167,0.289)
\or\qbezier(0.333,0)(-0.167,0.289)(-0.167,0.289)
\or\qbezier(-0.167,0.289)(-0.167,0.289)(-0.167,0.289)
\or\qbezier(-0.167,-0.289)(-0.167,0.289)(-0.167,0.289)
\or\qbezier(0,0)(-0.167,0.289)(-0.167,0.289)
\fi
\or\ifcase#2
\qbezier(0.667,0)(-0.167,-0.289)(-0.167,-0.289)
\or\qbezier(0.206,0.634)(-0.167,-0.289)(-0.167,-0.289)
\or\qbezier(-0.539,0.392)(-0.167,-0.289)(-0.167,-0.289)
\or\qbezier(-0.539,-0.391)(-0.167,-0.289)(-0.167,-0.289)
\or\qbezier(0.206,-0.634)(-0.167,-0.289)(-0.167,-0.289)
\or\qbezier(0.333,0)(-0.167,-0.289)(-0.167,-0.289)
\or\qbezier(-0.167,0.289)(-0.167,-0.289)(-0.167,-0.289)
\or\qbezier(-0.167,-0.289)(-0.167,-0.289)(-0.167,-0.289)
\or\qbezier(0,0)(-0.167,-0.289)(-0.167,-0.289)
\fi
\or\ifcase#2
\qbezier(0.667,0)(0,0)(0,0)
\or\qbezier(0.206,0.634)(0,0)(0,0)
\or\qbezier(-0.539,0.392)(0,0)(0,0)
\or\qbezier(-0.539,-0.391)(0,0)(0,0)
\or\qbezier(0.206,-0.634)(0,0)(0,0)
\or\qbezier(0.333,0)(0,0)(0,0)
\or\qbezier(-0.167,0.289)(0,0)(0,0)
\or\qbezier(-0.167,-0.289)(0,0)(0,0)
\or\qbezier(0,0)(0,0)(0,0)
\fi
\fi
}

\def\Trythis{
\qbezier(-1,0)(-1,0)(0.866,-0.5)
\put(-0.0536,-0.2){\circle*{0.15}}
}


\def\DoDash[#1,#2]{
\ifcase#1
\ifcase#2
\qbezier[20](0.5,0)(0,-1.25)(-0.5,0)
\or\qbezier[20](0.5,0)(-0.25,-1.25)(-1,-0.5)
\or\qbezier[20](0.5,0)(+0.25,-1.25)(-0.25,-0.5)
\or\qbezier[20](0.5,0)(0,-2.0)(-1.5,-1)
\or\qbezier[20](0.5,0)(0,-2.0)(-0.75,-1)
\or\qbezier[20](0.5,0)(0.25,-2.0)(-0.5,-1)
\or\qbezier[20](0.5,0)(0.25,-2.0)(0,-1)
\fi
\or\ifcase#2
\qbezier[20](0.25,-0.5)(-0.25,-1.25)(-0.5,0)
\or\qbezier[20](0.25,-0.5)(-0.25,-1.25)(-1,-0.5)
\or\qbezier[20](0.25,-0.5)(0,-1.25)(-0.25,-0.5)
\or\qbezier[20](0.25,-0.5)(0,-2.0)(-1.5,-1)
\or\qbezier[20](0.25,-0.5)(0,-2.0)(-0.75,-1)
\or\qbezier[20](0.25,-0.5)(0,-2.0)(-0.5,-1)
\or\qbezier[20](0.25,-0.5)(0.25,-2.0)(0,-1)
\fi\or\ifcase#2
\qbezier[20](1,-0.5)(0,-1.5)(-0.5,0)
\or\qbezier[20](1,-0.5)(0,-1.5)(-1,-0.5)
\or\qbezier[20](1,-0.5)(0,-1.5)(-0.25,-0.5)
\or\qbezier[20](1,-0.5)(0,-2.0)(-1.5,-1)
\or\qbezier[20](1,-0.5)(0,-2.0)(-0.75,-1)
\or\qbezier[20](1,-0.5)(0,-2.0)(-0.5,-1)
\or\qbezier[20](1,-0.5)(0,-2.0)(0,-1)
\fi
\fi
}

\def\DoA{
\put(0,1){\line(0,-1){0.5}}
}
\def\DoB{
\put(-0.5,0){\line(1,1){0.5}}
}\
\def\DoC{
\put(0,0.5){\line(1,-1){0.5}}
}
\def\DoD{
\put(-1,-0.5){\line(1,1){0.5}}
}
\def\DoE{
\put(-0.5,0){\line(1,-2){0.25}}
}
\def\DoF{
\put(0.25,-0.5){\line(1,2){0.25}}
}
\def\DoG{
\put(0.5,0){\line(1,-1){0.5}}
}
\def\DoH{
\put(-1.5,-1){\line(1,1){0.5}}
}
\def\DoI{
\put(-1,-0.5){\line(1,-2){0.25}}
}
\def\DoJ{
\put(-0.5,-1){\line(1,2){0.25}}
}
\def\DoK{
\put(-0.25,-0.5){\line(1,-2){0.25}}
}
\newtheorem{defn}{Definition}[section]
\newtheorem{theorem}[defn]{Theorem}
\newtheorem{lemma}[defn]{Lemma}
\newtheorem{corollary}[defn]{Corollary}
\newtheorem{fact}[defn]{Fact}
\newtheorem{prop}[defn]{Proposition}
\newtheorem{example}[defn]{Example}
\newtheorem{conj}[defn]{Conjecture}
\newtheorem{question}[defn]{Question}

\newcommand{\rtb}{\begin{flushright}
\begin{picture}(0.3,0.3)(0,0)
\put(-0.3,0){\framebox(0.3,0.3)}
\end{picture}
\end{flushright}}
\newcommand{\hb}{\hbar}
\newcommand{\Ch}{ {\cal A} }
\newcommand{\ChF}{\hat{\Ch}}
\newcommand{\bq}{\begin{equation}}
\newcommand{\ci}{\epsilon}
\newcommand{\eq}{\end{equation}}
\newcommand{\ra}{\rightarrow}
\newcommand{\ls}{sl(2,{\cal C})}
\newcommand{\oc}{\overline{c}}
\newcommand{\ve}{\varepsilon}
\newcommand{\gl}{gl(1|1)}
\newcommand{\ov}{\overline{1}}
\newcommand{\zv}{\overline{0}}
\newcommand{\EJ}{\hat{J}_{sl(2),\lambda}}

\begin{titlepage}
\hskip 1cm
\begin{flushright}
q-alg/9603011\\
\vspace{12pt}
February 1996
\end{flushright}

\hskip 1.5cm
\vfill
\begin{center}
{\LARGE Alexander-Conway limits of many Vassiliev weight systems.
}\vfill
{\large  
A.Kricker
\vspace{6pt}
}\\[10mm]
{\it School of Mathematical Sciences, \\
 University of Melbourne, Parkville, 3052, \\
Australia.}
\end{center}
\vfill
\begin{center}
{\bf Abstract}
\end{center}
{\small  
Previous work has shown that certain leading orders of arbitrary Vassiliev
invariants are generically in the algebra of the coefficients of the
Alexander-Conway polynomial \cite{KSA}. Here we illustrate this for a large
class of examples, exposing the simple logic behind several existing results
in the literature \cite{FKV,BNG}. This approach facilitates an extension to a
large class of Lie (super)algebras.
}
\vfill \hrule width 3.cm
{\footnotesize
\noindent
$^1$
Email: krick@mozart.ph.unimelb.edu.au}

\end{titlepage}

\section{Introduction.}

The study of Vassiliev knot invariants is something of a hybrid enterprise.
Defined by a topologist  in 1990 \cite{Vas} it was soon realised  that they
provided common ground for knot invariants arising from disparate
investigations in statistical mechanics, quantum groups, topological quantum
field theories and algebraic topology \cite{BN}.

The characterising feature of the Vassiliev infrastructure is the equivalence
of many topological and algebraic statements regarding Vassiliev knot
invariants to certain assertions defined in terms of combinatorics and
graph theory. This relationship was successfully exploited in \cite{BNG} to
provide a proof of a conjecture due to Melvin, Morton and Rozansky.

 We denote by $\EJ(\hb)[.]$ 
the framing independent 
$U_q(sl(2))$ invariant evaluated in a representation of dimension $\lambda
+ 1$ at $q\ =\ e^{\hb}$. Write the Alexander-Conway  polynomial,
$C(z)[.]$ (see section \ref{AC}).

\begin{theorem}[\cite{BNG}]
Expanding $\EJ/(\lambda + 1)$ in powers of $\lambda$ and $\hb$,
\bq
\frac{\EJ(K)(\hb)}{\lambda + 1}\ =\ \sum_{j,m\geq 0} b_{jm}(K)\lambda^j \hb^m,
\eq
we have,
\begin{enumerate}
\item{$b_{jm}(K) = 0$ if $j>m$.}
\item{Define \[ JJ(K)(\hb) = \sum_{m=0}^\infty b_{mm}(K) \hb^m. \]

        Then,
\bq
JJ(K)(\hb).\frac{\hb}{e^{\frac{\hb}{2}}-e^{-\frac{\hb}{2}}}C(K)(e^\frac{\hb}{2}-e^\frac{-\hbar}{2}) = 1.
\eq
}
\end{enumerate}
\rtb
\end{theorem}

In \cite{BNG} there is proved a generalisation of the MMR conjecture to the quantum group
invariants arising from general semi-simple Lie algebras: the series of
 terms of highest order in $\lambda$ turns out to be a product of inverse 
Alexander-Conway polynomials, one for each positive root, with weighted
indeterminates.  

Following a further suggestion in \cite{BNG}, \cite{KSA} demonstrated the universality
of this result: that arbitrary Vassiliev invariants had a ``highest order
piece'' which was always constructible from the coefficients of the
Alexander-Conway polynomial. Moreover, \cite{KSA} showed that this piece was
naturally selected by the action of the Adams operations (cabling operations).

In another investigation, \cite{FKV} uncovered the sequence of
Alexander-Conway weight systems in the universal weight system of a Lie
superalgebra, $gl(1|1)$. $W_{\gl} \subset Z(U(\gl))$ which is a polynomial
ring on two commuting generators, $c$ and $h$. They discovered that if one
``set by hand'' $c$ to 0 and $h$ to 1, then $W_{\gl}$ reproduced the sequence
of Alexander-Conway weight systems. This was achieved via the representation
theory of the cyclic $\gl$-modules, which was used to derive a recursive
relation for the evaluation of $W_{\gl}$. An argument that the $c\mapsto 0$
prescription is equivalent to deframing appears in \cite{FKV}, and was
indicated in \cite{KSA}.

In this work, we use the results of \cite{KSA} to expose the underlying
simplicity (and unity) of the results of \cite{BNG} and
\cite{FKV}. Essentially, we employ a description of a basis for the chord
diagram algebra which isolates the subspace on which the sequence of
Alexander-Conway weight systems are non-zero. This basis was introduced in 
\cite{KSA}.

In particular:
\begin{enumerate}
\item{We provide a simple proof of the Melvin-Morton-Rozansky
conjecture.}
\item{ We see that $\hat{W}_{\gl}$ is precisely the sequence of
Alexander-Conway weight systems, the recipe $c\mapsto 0$ being implicit in the
requirement of framing independence.  }
\item{ We combine elements of the above two results in the proof of the
following, new result. This result contains the generalisation of MMR to
semi-simple Lie algberas of \cite{BNG}.
\begin{theorem}
Take $L=L_{\zv}\oplus L_{\ov}$, a Lie superalgebra of classical type with
semi-simple body. With regard to a Cartan subalgbera $H\subset L_{\zv}$, $L$
has a root decomposition $L=\oplus_{\alpha\in \Delta} L_\alpha$. The set of
roots $\Delta$ partitions $\Delta = \Delta_{\zv}\cup \Delta_{\ov}$. Regarding
$\hat{W}_{L,\lambda,n}:\Ch_n\ra {\cal C}$ as a weight system valued function of the
highest weight $\{\lambda\}$ of some irreducible $L$-module:
\begin{enumerate}
\item{$\frac{\hat{W}_{L,\lambda,n}}{dim V_{\lambda}}$ is polynomial in $\{\lambda\}$
of highest order $n$. Write the coefficient of the order $n$ term $k_{nn}$.}
\item{
\bq
\{ \sum_n k_{nn} \hbar^n \} \circ Z_K\ =\ \frac{\prod_{\alpha\in \Delta^+_1}
\tilde{C}(<\lambda,\alpha>\hbar)}{\prod_{\beta\in \Delta^+_0}
\tilde{C}(<\lambda,\beta>\hbar )},
\eq
}
\end{enumerate}
where in the above $\tilde{C}$ is the normalised Alexander-Conway polynomial
(section \ref{AC}), $Z_K$ is the Kontsevich functor (section \ref{AC}), and
the precise meaning of the above statements is discussed in section
\ref{super}. 
\rtb
\end{theorem}
}
\end{enumerate}

Can we identify this as the limit of some ``quantum'' knot invariant?
Knot invariants can be defined from quantum supergroups
$V_{U_q(L)}$ \cite{susy}. 
\begin{question}
Is $V_{U_q(L)}$ a canonical knot invariant, with weight system that
constructed from the Lie superalgebra $L$ (as in \cite{Vain}, and reviewed
in section \ref{weights}) ?
\rtb
\end{question}
This requires some (perhaps trivial) extension of Kohno's work \cite{Koh}.

This work is timely for several reasons. There has recently been some interest
displayed in investigating the sub-leading terms of the MMR series. In
particular, Rozansky has formulated several intriguing conjectures
\cite{R2,Gar}.
 As the
MMR series is canonical, these lower order terms should be most easily
attacked through the formalism of weight systems. A proper understanding of
the leading order terms is necessary before this extension can be attempted.

In other recent work \cite{Vog}, Vogel
 has demonstrated that the space of weight systems
originating in Lie superalgebras is strictly larger than the space of weight
systems coming from semi-simple Lie algebras. The work presented here
 is interesting as 
an indication of the relationship between the two spaces.

In section 1 we exhibit the theory relevant for our analyses. In particular,
the canonical deframing operation and weight systems from Lie superalgebras
are recalled. In section 2 we examine several cases: the
MMR conjecture, the $gl(1|1)$ weight system and the class of Lie superalgebras
of classical type. 

\section{Review of results.}

\subsection{Notation.}

In all that follows, we shall work over the complex numbers ${\cal C}$.
Denote by ${\cal Z}_+$ the set of positive integers including zero.

Recall that a Chinese character diagram (CCD)
is a finite graph with trivalent and univalent
vertices such that the univalent vertices are cyclically ordered on a copy
of $S^1$ (the Wilson loop) and such that the incoming
edges to a trivalent vertex are cyclically ordered.
 We denote $\ChF_*$ for the graded
 ${\cal C}$-vector space
 generated by Chinese character diagrams. A chord diagram is
graded by half the number of vertices in it's associated graph (we shall refer
to a chord diagram with 2n such vertices as an $n-$CCD). 
Recall that the STU subspace of $\ChF_n$ is generated by expressions
on $n-$CCDs of the following form:
\bq
\Picture{\Sterm}\ -\ \Picture{\Tterm}\ +\ \Picture{\Uterm}\ .
\eq

\

We denote $\Ch_n$ for the graded ${\cal C}$-vector space that results from
the quotient of $\ChF_n$ by the STU subspace on $n-$CCDs. This is naturally
a Hopf algebra with commutative product, and co-commutative co-product. The
product occurs via connect sum, and respects the grading. The coproduct is 
as follows. Denote the set of components of the associated graph of an
$n-$CCD $X$ via $C(X)$. Denote the CCD obtained from $X$ by deleting some
of the connected components  $\{ a_1,a_2,\ldots\}$ by $X_{\{a_1,a_2,\ldots\}}$
. Then
\bq
\Delta(X)\ =\ \sum_{S\subset C(X)} X_S \otimes X_{C(X)\backslash S}.
\eq

An $n-$weight system is a linear map $\Ch_n \ra {\cal C}$. Knot invariants
of finite order $n$ give well-defined $n-$weight systems. This is
achieved by evaluating on an $n-$knot representative of a given chord diagram,
in the usual fashion \cite{BN}. We multiply $X$ and $Y$, weight systems of
order $n$ and $m$ respectively, by $X.Y(v)= X\times Y(\Delta(v)),\ v\in
\Ch_{n+m}$. 

There exists an alternative basis for $\Ch_*$ which is often very useful:
these are the Chinese characters (CC). Chinese characters are graphs 
with univalent and trivalent vertices. They have no Wilson loop and
are again graded by half the number of vertices. An $n-$CC $v$ represents the sum
with unit coefficients over ways of locating the univalent vertices on a
Wilson loop. We will refer to the set of summands in the above sum as $E(v)$
(i.e. this is a set of $n-$CCDs). 
The STU relations in $\Ch_*$ pull back to the IHX and antisymmetry
relations of CCs. 

{\bf Antisymmetry.}
\begin{equation}
\Picture{
\put(0,0){\vector(0,1){1}}
\put(0,0){\vector(1,-1){1}}
\put(0,0){\vector(-1,-1){1}}
}\ \ =\ \
\ \ -\ \Picture{
\put(0,0){\vector(0,1){1}}
\qbezier(0,0)(-1,-0.5)(1,-1)
\qbezier(0,0)(1,-0.5)(-1,-1)
\put(1,-1){\vector(2,-1){0.1}}
\put(-1,-1){\vector(-2,-1){0.1}}
}
\end{equation}

\

{\bf IHX relation.}
\begin{equation}
\Picture{
\put(-1,1){\vector(1,0){2}}
\put(-1,-1){\vector(1,0){2}}
\put(-1,1){\vector(-1,0){0.1}}
\put(-1,-1){\vector(-1,0){0.1}}
\put(0,1){\line(0,-1){2}}
}
\ \ =\ \ 
\Picture{
\qbezier(-1,1)(-0.5,1)(-0.5,0)
\qbezier(-1,-1)(-0.5,-1)(-0.5,0)
\qbezier(1,1)(0.5,1)(0.5,0)
\qbezier(1,-1)(0.5,-1)(0.5,0)
\put(-0.5,0){\line(1,0){1}}
\put(-1,1){\vector(-1,0){0.1}}
\put(-1,-1){\vector(-1,0){0.1}}
\put(1,1){\vector(1,0){0.1}}
\put(1,-1){\vector(1,0){0.1}}
}\ \ -\ \
\Picture{
\qbezier(-1,1)(-0.5,1)(-0.25,0)
\qbezier(-0.25,0)(0,-1)(1,-1)
\qbezier(1,1)(0.5,1)(0.25,0)
\qbezier(-1,-1)(0,-1)(0.25,0)
\put(-0.25,0){\line(1,0){0.5}}
\put(-1,1){\vector(-1,0){0.01}}
\put(-1,-1){\vector(-1,0){0.01}}
\put(-1,-1){\vector(-1,0){0.01}}
\put(1,1){\vector(1,0){0.01}}
\put(1,-1){\vector(1,0){0.01}}
}
\end{equation}

\

\begin{lemma}[\cite{BN}]
The subspace generated by $n-$CCs spans the space of $n-$CCDs in $\Ch_n$.
\rtb
\end{lemma}

\subsection{The deframing projector.}
The weight systems produced by knot invariants
 satisfy an additional relation,
the one-term (or 1-T) relation \cite{BN}. Writing the $n-$weight system
obtained from a finite type knot invariant $V$ of order $n$ $W_n[V]$, then:

\begin{equation}
W[V] \left(
\Picture{
\DottedCircle
\FullChord[6,7]
\Arc[6]
\Arc[7]}
\right)
\ =\  
V\left(
\Picture{
\qbezier(-0.8,0.9)(-0.8,0.45)(0,0)
\qbezier(-0.8,-0.9)(-0.8,-0.45)(0,0)
\put(0,0){\circle*{0.15}}
\qbezier(0,0)(1,0.5)(1,0)
\qbezier(0,0)(1,-0.5)(1,0)}
\right)
\  =\ 
V\left(
\Picture{
\qbezier(-0.8,0.9)(-0.8,0.45)(-0.2,0.1)
\qbezier(0.2,-0.10)(1,-0.5)(1,0)
\qbezier(-0.8,-0.9)(-0.8,-0.45)(0,0)
\qbezier(0,0)(1,0.5)(1,0)
}\right)\ -\
V\left(
\Picture{
\qbezier(-0.8,0.9)(-0.8,0.45)(0,0)
\qbezier(0,0)(1,-0.5)(1,0)
\qbezier(-0.8,-0.9)(-0.8,-0.45)(-0.2,-0.1)
\qbezier(0.2,0.1)(1,0.5)(1,0)
}\right)\  =\ 0.\
\end{equation}

From an arbitrary weight system $W_n$ we may wish to construct a
 weight system most
representative of the original, but satisfying in addition the 1-T relation.
We do this by composing $W_n$ with the {\it deframing projector}, $\phi_n$
 \cite{BN}. To define $\phi_n$ we first introduce $s_n:\Ch_n\ra \Ch_{n-1}$.

\begin{defn}
Let $X\in \Ch_n$. 
\bq
s_n(X)\ =\ \sum_{chords\,  a \in C(X) } X_a\ \in \Ch_{n-1}.
\eq
A component of the graph of X is a chord if it has no trivalent vertices.
\rtb
\end{defn}

\begin{defn}
The {\it deframing projector} $\phi_n:\Ch_n \ra \Ch_n$ is constructed:
\bq
\phi_n=\ Id\ -\ \theta.s\ +\ \frac{\theta^2.s^2}{2!} -\ \ldots\ +\frac{(-1)^n
\theta^n.s^n}{n!}.
\eq
In the above, $\theta$ represents the diagram with a single chord, and the
product is the natural connect sum product in the graded chord diagram
algebra, $\Ch_*$. By $s^i$ we mean $s_{n-i+1}\circ\ldots\circ s_{n-1}\circ
s_n$. 
\rtb
\end{defn}

\begin{lemma}[\cite{BN,KSA}]
The operator $\phi_n$ satisfies the following properties:
\begin{enumerate}
\item{ $\phi_n$ is well-defined as a linear map $\phi_n:\Ch_n \ra \Ch_n$,}
\item{ $\phi_n(
\Picture{ \Arc[5]\Arc[6]\DottedCircle\FullChord[5,6]} )\ =\ 0,$ }
\item{ $\phi_n(v)\ =\ v\ $ if and only if $s_n(v)\ =\ 0$, for $v\in \Ch_n$.}
\end{enumerate}
\label{DIS}
\rtb
\end{lemma}

Denote by  $I_n$ the deframing invariant subspace of $\Ch_n\ $ ($I_n=\{v\in
\Ch_n: \phi_n(v)=v\}$).
\begin{lemma}[\cite{KSA}]
\bq
I_n\ =\ Ker(s_n)\ =\ Span\{\mbox{$n-$CCs without chords.}\}
\eq
\rtb
\end{lemma}

The deframing invariant subspace is thus spanned by those Chinese characters 
which have trivalent vertices in every component. This collection of
generators
is partitioned by the number of univalent vertices a generator has. Write
$I_n^p$ for the subspace of $I_n$ generated by $n-$CCs without chords
with $p$ univalent
vertices.

\begin{lemma}[\cite{KSA}]
\bq
I_n\ =\ \oplus_{i=1}^n I^i_n.
\eq
\rtb
\end{lemma}

There exists a description of a basis for the ``highest'' piece $I_n^n$.

\begin{defn}
Select a partition of $n=n_1 + n_2 + \ldots + n_j, n_i\in {\cal Z}^+$. 
This is represented by the $j$-tuplet $\{n_1,\ldots,n_j\}$.  Each summand
 in this
partition corresponds to a component in the $n-$CC $\tau_{\{n_1,\ldots,n_j\}}$.
 The
component  corresponding to $n_i$ is a loop of $n_i$ edges with extra edges
attached radially at each vertex. For example:
\[
\tau_{[2]}\ =\
\Picture{
\put(0,1){\line(0,-1){0.6}}
\put(0,-0.4){\line(0,-1){0.6}}
\qbezier(0,0.4)(-0.7,0)(0,-0.4)
\qbezier(0,0.4)(0.7,0)(0,-0.4)
}\ \ ,\ \
\tau_{[4,2]}\ =\
\Picture{
\put(-1,1){\line(0,-1){2}}
\put(0,1){\line(0,-1){2}}
\put(-1,0.4){\line(1,0){1}}
\put(-1,-0.4){\line(1,0){1}}
\put(1,1){\line(0,-1){0.6}}
\put(1,-0.4){\line(0,-1){0.6}}
\qbezier(1,0.4)(0,0)(1,-0.4)
\qbezier(1,0.4)(2,0)(1,-0.4)
}\ \ \ \ .\]
\rtb
\end{defn}

\begin{lemma}[\cite{KSA}]
\bq
I^n_n\ =\ \mbox{span}\{ \{ \tau_{\bf P} \}\ ;\ {\bf
P}=\{p_1,\ldots,p_j\}:p_1+\ldots+p_j=n\ \&\ p_i\in 2{\cal Z}^+  \}
\eq
These ``even partitions'' form a linearly independent basis for $I^n_n$.
\rtb
\end{lemma}  

For example:
\begin{eqnarray}
I^2_2 & = & \mbox{span}\{ 
\Picture{
\put(0,1){\line(0,-1){0.6}}
\put(0,-0.4){\line(0,-1){0.6}}
\qbezier(0,0.4)(-0.7,0)(0,-0.4)
\qbezier(0,0.4)(0.7,0)(0,-0.4)
} \}\ ,\ \nonumber \\
& & \nonumber \\
I^4_4 & = & \mbox{span}\{
\Picture{
\put(-0.5,1){\line(0,-1){0.6}}
\put(-0.5,-0.4){\line(0,-1){0.6}}
\qbezier(-0.5,0.4)(-1,0)(-0.5,-0.4)
\qbezier(-0.5,0.4)(0,0)(-0.5,-0.4)
\put(0.5,1){\line(0,-1){0.6}}
\put(0.5,-0.4){\line(0,-1){0.6}}
\qbezier(0.5,0.4)(0,0)(0.5,-0.4)
\qbezier(0.5,0.4)(1,0)(0.5,-0.4)
}\ \ ,\ \
\Picture{
\put(-0.5,1){\line(0,-1){2}}
\put(0.5,1){\line(0,-1){2}}
\put(-0.5,0.4){\line(1,0){1}}
\put(-0.5,-0.4){\line(1,0){1}}
} \}\ ,\  \nonumber \\
& & \nonumber \\
I^6_6 & = & \mbox{span}\{ \ 
\Picture{
\put(-1,1){\line(0,-1){0.6}}
\put(-1,-0.4){\line(0,-1){0.6}}
\qbezier(-1,0.4)(-1.3,0)(-1,-0.4)
\qbezier(-1,0.4)(-0.7,0)(-1,-0.4)
\put(0,1){\line(0,-1){0.6}}
\put(0,-0.4){\line(0,-1){0.6}}
\qbezier(0,0.4)(0.3,0)(0,-0.4)
\qbezier(0,0.4)(-0.3,0)(0,-0.4)
\put(1,1){\line(0,-1){0.6}}
\put(1,-0.4){\line(0,-1){0.6}}
\qbezier(1,0.4)(1.3,0)(1,-0.4)
\qbezier(1,0.4)(0.7,0)(1,-0.4)
}\ \ ,\ \
\Picture{
\put(-1,1){\line(0,-1){2}}
\put(0,1){\line(0,-1){2}}
\put(-1,0.4){\line(1,0){1}}
\put(-1,-0.4){\line(1,0){1}}
\put(1,1){\line(0,-1){0.6}}
\put(1,-0.4){\line(0,-1){0.6}}
\qbezier(1,0.4)(0,0)(1,-0.4)
\qbezier(1,0.4)(2,0)(1,-0.4)
}\ \ \ ,\ \
\Picture{
\put(-1,1){\line(0,-1){2}}
\put(1,1){\line(0,-1){2}}
\put(0,-0.4){\line(0,-1){0.6}}
\put(0,1){\line(0,-1){0.6}}
\put(-1,0.4){\line(1,0){2}}
\put(-1,-0.4){\line(1,0){2}}
}\  \}\ \ .\ 
\end{eqnarray}

\

In \cite{KSA} the set of $n-$CCDs $E(\tau_{\bf P})$ were referred to as
``CCDs with all graph components basic loops''. $I^*_*\,(\equiv \oplus_i I^i_i)$ is a graded,
commutative and cocommutative Hopf algebra, with the product being disjoint
union and the usual coproduct. Note that \cite{KSA} essentially showed that
$I^*_*$ is the dual Hopf algbera to the Hopf algebra of immanent weight
systems \cite{BNG}.

\subsection{Alexander-Conway type weight systems.}
\label{AC}

\begin{defn}
The Alexander-Conway polynomial $C(z)[.]\in {\cal Z}[z]$ is a
knot invariant defined by:
\begin{enumerate}
\item{ the skein relation 
\bq
C(z)\left( \Picture{
\DottedCircle
\put(-0.707,0.707){\vector(1,-1){1.41}}
\put(-0.707,-0.707){\line(1,1){0.6}}
\put(0.1,0.1){\vector(1,1){0.6}}
} \right)\ - \
C(z)\left( \Picture{
\DottedCircle
\put(-0.707,0.707){\line(1,-1){0.6}}
\put(-0.707,-0.707){\line(1,1){1.41}}
\put(0.1,-0.1){\vector(1,-1){0.6}}
\put(0.1,0.1){\vector(1,1){0.6}}
} \right)\ =\
z\ .\ C(z)\left( \Picture{
\DottedCircle
\qbezier(-0.707,0.707)(0,0.4)(0.707,0.707)
\qbezier(-0.707,-0.707)(0,-0.4)(0.707,-0.707)
\put(0.707,0.707){\vector(2,1){0.1}}
\put(0.707,-0.707){\vector(2,-1){0.1}}
} \right),
\eq
and}
\item{ 
the boundary condition
\bq C(z)\left(
\Picture{
\FullCircle
}
\right)=1.
\eq}
\end{enumerate}
\rtb
\end{defn}

Write $C(z)[.]=\sum \tilde{c}_n[.]z^n$. $\tilde{c}_n$ is a knot invariant
of finite order $n$. As such it canonically gives rise to an $n-$weight system
-- $W_n(\tilde{c}_n)$ -- which we shall write $c_n:\Ch_n\ra {\cal C}$. Through
their theory of intersection graphs and cycle decomposition sums \cite{BNG}
provided a powerful 
characterisation of $c_n$. In \cite{KSA} the
authors used that characterisation to produce the following realisation
of the $\{c_n\}$ on the $I_n^p$:
\begin{theorem}
\label{cnvalues}
Take $v\in \Ch_n$.
\bq
c_n(v)\ =\ \left\{ \begin{array}{ll}
0 & \mbox{if}\ \phi_n(v)\ =\ 0, \\
0 & \mbox{if}\ v\ \in\ I^{<n}_{n}, \\
(-2)^{\#({\bf P})} & \mbox{if}\ v\ \in\ E(\tau_{\bf P}).
\end{array}
\right.
\eq
In the above $\#(\tau_{\bf P})$ represents the number of terms in the
partition ${\bf P}$. Note that the trivial case of an empty graph is
implicit in the above, $c_0$ taking the value $1$.
\rtb
\end{theorem}

 We may ask whether any weight system vanishing on $Ker\,\phi_n$ and on
$I^{<n}_{n}$ is in 
the algebra of Alexander-Conway weight systems (a sum of products of the
$c_n$, hereafter referred to as the c-algebra)? The affirmative answer is
implicit from \cite{BNG} and 
\cite{KSA}. For the present approach it is relevant that we include a short
argument of this fact.

\begin{lemma}
Any $n-$weight system vanishing on Ker\,$\phi_n$ and on $I_n^{<n}$ is equal to
a sum
of products of the $c_n$. 
\label{cnalg}
\

{\bf \underline{Proof.}}
A few points to start. Write $W_{\bf P}$ for the weight system vanishing on
$Ker\,\phi_n$ and on $I^{<n}_n$, and
orthonormal to $\tau_{\bf P}$ on $I^n_n$. We have that $W_{\{n_1,n_2,\ldots,n_j\}}=
W_{\{n_1\}}.W_{\{n_2\}}\ldots W_{\{n_j\}}.$ (These span the Hopf algebra
dual to $I^*_*$). 

The proof proceeds by induction on n, the grade of $\Ch_n$. Assume it is true
that the dual basis for $I^m_m,m<n$ is in the c-algebra. Then it
is true that for all $\tau_{\bf P}\neq \tau_{ \{ n \} }$ in $I^n_n$ 
the dual weight systems $W_{\bf P}$ are
products
of weight systems of lower orders, hence by assumption in the c-algebra.
  It remains to show that $W_{ \{n\} }$ is in the c-algebra. $c_n$
is non-vanishing on $\tau_{\{n\}}$. Thus $W_{\{n\}}=kc_n+\sum_{{\bf P}\neq
\{n\}} W_{\bf P}k_{\bf P}$ for some constants $k$. This completes the
induction.  It is trivial for $n=0$ so the proof is complete.
\rtb
\end{lemma}

Kontsevich has defined a ``universal knot invariant'' taking values in the
algebra of chord diagrams, $Z_K:\{\mbox{Knots}\}\ra \Ch_*$. It is universal
 in the sense that given any knot invariant invariant $V$ of finite order
$n$, we have:
\bq
V\ =\ W_n[V]\circ Z_K\ \ +\ \ \left( 
\begin{array}{l} \mbox{knot invariants of} \\
\mbox{finite order}\ <n.
\end{array}
\right)\ .
\eq
In the case that $V=W_n[V]\circ Z_K$ we call $V$ a canonical Vassiliev
invariant.
Note that $\Delta(Z_K)=Z_K\otimes Z_K$. Thus the product of two canonical
Vassiliev invariants is canonical with multiplied weight systems.
 Note further that $C(\hbar)\neq c(\hbar)\circ Z_K (= c_i\hbar^i\circ Z_K).$ \cite{BNG}
found a normalised form of the Alexander-Conway polynomial which {\it is}
 canonical.
\begin{prop}[\cite{BNG}]
\bq
c(\hbar)\circ Z_K\ =\ \frac{\hbar}{(e^{\frac{\hbar}{2}}\ -\
e^\frac{-\hbar}{2})}C(e^{\frac{\hbar}{2}}\ -\ e^{-\frac{\hbar}{2}})\ \ \equiv\
\tilde{C}(\hbar). 
\eq
\rtb
\end{prop}

It makes sense to talk about the inverse of the series of weight systems $c$. We must
first introduce some notion of identity. The identity knot invariant is
canonical of order 0.
\begin{defn}
Define a sequence of weight systems $\ve_*:\Ch_*\ra {\cal C}$ by the
following. Take $v\in \Ch_n$:
\bq
\ve_n (v)\ =\ \left\{ 
\begin{array}{ll}
1 & \mbox{if}\ n=0, \\
0 & \mbox{otherwise}.
\end{array} \right.
\eq
\rtb
\end{defn}

\begin{defn}
\label{inv}
Define a sequence of weight systems $\oc_*:\Ch_*\ra{\cal C}$ by the
following. Take $v\in \Ch_n:$
\bq
\oc_n(v)\ =\ \left\{ \begin{array}{ll}
0 & \mbox{if}\ \phi_n(v)\ =\ 0, \\
0 & \mbox{if}\ v\ \in\ I^{<n}_{n}, \\
(2)^{\#({\bf P})} & \mbox{if}\ v\ \in\ E(\tau_{\bf P}).
\end{array}
\right.
\eq
\rtb
\end{defn}

\begin{lemma}
\bq
\ve\ =\ \sum_{i=0}^{n}c_i.\oc_{n-i}\ =\ \sum_{i=0}^{n}\oc_{n-i}.c_i.
\eq
{\bf \underline{Proof.}}
The second equality is immediate from the commutativity of the coproduct.
The case $n=0$ is immediate. What remains to be shown is
 that $\sum_{i=0}^{n}(c_i\times \oc_{n-i})\Delta(v)\ =\
0$ for $v\in \Ch_{\geq 1}$.
$\Delta(v)$ is a sum of terms like $v_{X}\otimes v_{C(v)/X}$. If $v$ has 
j components and $X$ has order $x$ then $v_{X}$ has $j-x$ components and
$v_{C(v)/X}$ has $x$. This will be non-zero for the $c_{j-x}\times \oc_{x}$
term
which evaluates to $(-2)^{j-x}(2)^x$. So the whole expression evaluates to
$\sum_{i=0}^{j}\left( \begin{array}{l} j \\ i \end{array}
\right)(-2)^i(2)^{j-i}\ =\ (2-2)^j\ =\ 0.$
\rtb
\end{lemma}

The motivation for calling the sequence $\oc_*$ an inverse of $c_*$ is that
one can (formally) collect sequences $c(\hbar)=c_i \hbar^i$ and $\oc(\hbar) = \oc_i \hbar^i$
and write,
\bq
c.\oc\ =\ \ve.
\eq
Thus,
\bq
\oc \circ Z_K\ =\ \frac{1}{\tilde{C}(\hbar)}.
\eq

\subsection{Weight systems from Lie (super)algebras.}
\label{weights}
Here we recall how to construct sequences of weight systems (following
\cite{BN,Vain}) from familiar algebras. For example to build an $n-$weight
system we require a linear  mapping  $\ChF_n \ra {\cal C}$ which respects STU
relations.

In \cite{Vain}, a categorical formalism of
some generality is put forward for constructing weight systems:
 the theory of self-dual Lie S-algebras. We will expose here some elements
of that approach.

We recall first the definition of a Lie (super)algebra.
\begin{defn}
A Lie superalgebra is a (here, finite dimensional) ${\cal Z}_2$-graded vector
space $L=L_{\overline 0}\oplus L_{\overline 1}$ with a
product satisfying conditions given below. Write $[v]$ for the grade
of $v\in L$. Take $v,w,x\in L$. 
\begin{enumerate}
\item{$[L_i.L_j]\ \subset\ L_{i+j(mod 2)},$}
\item{$[v.w]\ =\ -(-1)^{[v][w]}[w.v],$}
\item{$[v.[w.x]]\ =\ [[v.w].x]\ +\ (-1)^{[v][w]}[w.[v.x]] $}
\end{enumerate}
\rtb
\end{defn}
In the case that $L_{\overline 1}=\emptyset$ this is a
Lie algebra. $L_{\zv}$ is called the {\it body} of the algebra, and
$L_{\ov}$ is called the {\it soul.}

Take some Lie (super)algebra $L$ admitting an invariant, supersymmetric,
 non-degenerate metric
$\kappa:L\otimes L\ra {\cal C}$ (recall that a supersymmetric metric
satisfies $\kappa(v,w)=(-1)^{[v][w]}\kappa(w,v))$.
The Cartan-Killing form ($Tr( Ad(x)Ad(y))$) is an example
of an invariant, (super)symmetric form. When the Cartan-Killing form of
a Lie algebra is non-degenerate, then that algebra is semi-simple.
If $\kappa$ is non-degenerate then
 it has an inverse. Label this $b:{\cal
C}\ra L\otimes L$. 
is semi-simple. There exists a class of Lie superalgebras
admitting invariant, supersymmetric, nondegenerate metrics:
 the superalgebras of classical type
\cite{Kac} (for example $sl(m|n),\ osp(m|n)$).

Our weight system will factor through the following algebra.
\begin{defn}
The {\it universal enveloping algebra} of $V$ -- $U(V)$ -- is the quotient
of the tensor algebra on $V$ by the ideal generated by expressions of the
form $x\otimes y -(-1)^{[x][y]}y \otimes x - [xy]$ for $x,y\in V$.
\rtb
\end{defn}

The first step in the construction is to model the set of CCDs in a suitable
category.

\begin{defn}
Denote {\bf FG} for the category of Feynman graphs. Obj$\{ {\bf FG} \}\ =\
{\cal Z}_+$. A morphism from $p$ to $q$ is a finite graph with
 univalent
and trivalent vertices  such that:
\begin{enumerate}
\item{there is a cyclic ordering to the edges meeting at each trivalent
vertex,}
\item{the set of univalent vertices partitions into two
ordered sets of order $p$ and $q$. Call the set of order $p$ the set of
incoming legs, and call the set of order $q$ the set of outgoing legs.}
\end{enumerate}
\rtb
\end{defn}

This category has a graphical interpretation. Morphisms from $p$ to $q$ are
are simply graph diagrams in the region $\{(x,y),0\leq y\leq 1\}$ such that
the set of incoming legs lie ordered with increasing $x$ on the
line $y=1$, and the set of outgoing legs lie similarly ordered on the line
 $y=0$ .

This is naturally a tensor category. As usual, for $A,B\in \mbox{Mor}\{ {\bf
FG} \},\ A\otimes B$ is represented by the diagram for $A$ with the diagram
for
$B$ placed to it's right. If Dom$(A)=\mbox{Ran}(B)$ then $A\circ B$ is constructed
by placing $B$ over $A$ (identifying appropriate univalent vertices) then 
contracting the whole diagram to fit in the region $0\leq y \leq 1$.

We {\it represent} this category in our (super)symmetric algebra. That is,
we assign vector spaces to objects, and linear transformations between them
to morphisms (Feynman graphs). Here, the object $n$ is sent to $L^{\otimes
n}$. Below we present a list of (irreducible) Feynman graphs that generate
{\bf FG} by tensor products and compositions. Once we constrain the
representation to preserve tensor products of objects and morphisms then
defining the transformations associated with these irreducible graphs 
determines the representation. 

Chinese character diagrams are distinguished morphisms in this category: the
set of morphisms with domain zero. To see this, break the Wilson loop of some
CCD at some point and lay it out as a linear chord diagram. Thus, we can take
a CCD, express it as an {\bf FG}-morphism, then it's representation in our
(super)symmetric algebra is a linear map ${\cal C}\ra U(V)$. Of course, it
takes more work to show it is well-defined on $\ChF_*$ and more still to
show that it descends to $\Ch_*$.

\begin{theorem}[\cite{Vain}]
The following list of morphisms generates {\bf FG}. Take some Lie
(super)algebra  $L$ with an invariant, (super)symmetric, non-degenerate
 inner product $\kappa$. Also presented is a well-defined representation of
${\bf FG}$ in the tensor algebra of $L$. Considered as a mapping from
$\ChF_*$, this representation descends to $\Ch_*$.

\begin{eqnarray*}
(a) & \Picture{
\qbezier[20](-1,1)(0,1)(1,1)
\qbezier[20](-1,-1)(0,-1)(1,-1)
\put(0,1){\line(0,-1){2}}} & {\bf Id}\ :\ V\ \ra\ V, \\
& & \\
& & \\
(b) &
\Picture{
\qbezier[20](-1,1)(0,1)(1,1)
\qbezier[20](-1,-1)(0,-1)(1,-1)
\qbezier(-0.5,1)(0,-0.2)(0.5,1)
} & {\bf \kappa}\ :\ V\otimes V\ \ra {\cal C}, \\
& & \\
& & \\
(c) &
\Picture{
\qbezier[20](-1,1)(0,1)(1,1)
\qbezier[20](-1,-1)(0,-1)(1,-1)
\qbezier(-0.5,-1)(0,-0.2)(0.5,-1)
} & {\bf b}\ :\ {\cal C}\ \ra\ V\otimes V, \\
& & \\
& & \\
(d) & \Picture{
\qbezier[20](-1,1)(0,1)(1,1)
\qbezier[20](-1,-1)(0,-1)(1,-1)
\put(0,-1){\line(0,1){1}}
\put(0,0){\line(1,2){0.5}}
\put(0,0){\line(-1,2){0.5}}
} & {\bf f}\ :\ V\otimes V\ \ra\ V,\ {\bf f}(v,w)\ =\ [v.w], \\
& & \\
& & \\
(e) & \Picture{
\qbezier[20](-1,1)(0,1)(1,1)
\qbezier[20](-1,-1)(0,-1)(1,-1)
\put(-0.5,-1){\line(1,2){1}}
\put(0.5,-1){\line(-1,2){1}}
} & {\bf X}\ :\ V\otimes V\ \ra\ V\otimes V,\ {\bf X}(v\otimes w)\ =\
(-1)^{[v][w]}w\otimes v.
\end{eqnarray*}
\rtb
\end{theorem}

Call this map $W_L:\Ch_*\ra U(L)$. \cite{Vain} shows that $Image(W_L)\subset
{\cal Z}(U(L))$. A weight system can then be recovered by
choosing an irreducible representation of $L$ and evaluating the trace of
the image of $W_L$. If the representation is indexed by $\Lambda$,
 then we shall
denote this weight system $W_{L,\Lambda}$. The above weight systems restricted
to $\Ch_n$ shall be written $W_{L,n}$ and $W_{L,\Lambda,n}$.

For calculational purposes it is often easiest to introduce some linearly
independent basis for $L$ and treat the weight system as a state sum. Assume
$L$ has such a basis $\{v_0,\ldots,v_n\}$ such that $\{v_0,\ldots,v_i\}\in
L_{\overline{0}}\ \&\ \{v_{i+1},\ldots,v_n\}\in L_{\overline{1}}$. Write
$[v_i.v_j]={\bf f}^k_{ij}v_k$; $\kappa(v_i,v_j)=\kappa_{ij}$; and ${\bf
b}={\bf b}^{ij}v_i
\otimes v_j$ (with Einstein summation). Take the graph of $D\in \Ch_n$ and add
vertices for every vertical turning point.  Then $W_L(D)$ is equal to a sum
over all ``colourings'' of the edges of this graph with each summand weighted
by an appropriate product of entries in the tensors
${\bf f}^k_{ij},\kappa_{ij},{\bf b}^{ij}$. Below we present an example of such a
calculation: 

\bq
W_L (\ \
\setlength{\unitlength}{40pt}
\Picture{
\put(-1,-1){\vector(1,0){2}}
\qbezier(-0.8,-1)(-0.8,0.5)(-0.4,0.5)
\qbezier(-0.4,0.5)(-0.2,0.5)(0,0)
\qbezier(0,0)(0.2,-0.5)(0,-0.7)
\qbezier(0.8,-1)(0.8,0.5)(0.4,0.5)
\qbezier(0.4,0.5)(0.2,0.5)(0,0)
\qbezier(0,0)(-0.2,-0.5)(0,-0.7)
\put(0,-1){\line(0,1){0.3}}
\put(0,-0.7){\circle*{0.15}}
\put(0,0){\circle*{0.15}}
\put(-0.4,0.5){\circle*{0.15}}
\put(0.4,0.5){\circle*{0.15}}
\put(-1,0){\mbox{a}}
\put(-0.3,0.1){\mbox{b}}
\put(0.2,0.1){\mbox{c}}
\put(0.8,0){\mbox{d}}
\put(0.25,-0.5){\mbox{f}}
\put(-0.35,-0.5){\mbox{e}}
\put(0.1,-0.9){\mbox{g}}
}\ \
\setlength{\unitlength}{20pt}
)\ =\ \sum_{a,\ldots,g}^{n} {\bf b}^{ab}{\bf b}^{cd}\delta^{f}_{b}
\delta^{e}_{c}(-1)^{[b][c]}{\bf f}^{g}_{ef}\  (v_a \otimes v_g \otimes v_d).
\eq

\

\section{Alexander-Conway limits.}

\subsection{The Melvin-Morton-Rozansky conjecture.}

Take the simple Lie algebra $sl(2,{\cal C})$. This is a vector space on three
generators $\{ H,X,Y \}$ such that:
\bq [H.X] =  2X,\ \  [H.Y]  =  -2Y,\ \ [X.Y]  =  H\ .
\eq

There exists an invariant, symmetric, non-degenerate inner product $\kappa$ on
$\ls$ (written $\kappa(v,w)=<v,w>$). The non-vanishing products on
 basis elements are the following:
\bq
\frac{1}{2}<H,H>\ =\ <X,Y>\ =\ <Y,X>\ =\ 1.
\eq
The inverse of this form is easily constructed and is:
\bq
{\bf b}\ =\ \frac{1}{2}H\otimes H\ +\ X\otimes Y\ +\ Y\otimes X.
\eq

The set of irreducible representations 
are indexed by $\lambda \in {\cal Z}^+\backslash \{0\}$, and have dimension
$\lambda+1$.  On some set of generators $\{ v_0,v_1,\ldots,v_\lambda\}$ the
representation $V_{\lambda}$ may be presented $\rho_{\lambda}
:\ls\ra Gl(V_\lambda)$,
\begin{eqnarray}
\rho_{\lambda}(h).v_j & = &(\lambda-2j)v_j, \nonumber \\
\rho_{\lambda}(y).v_j & = &(j+1)v_{j+1}, \nonumber \\
\rho_{\lambda}(x).v_j & = &(\lambda -j+1)v_{j-1}.
\label{rep}
\end{eqnarray}  

We shall be investigating the deframed weight system
$\hat{W}_{\ls,\lambda}$. It coincides with
$W_{\ls,\lambda}$  
on the deframing invariant subspace $I_*$. Our approach is to
consider $W_{\ls,\lambda}$ as a weight system valued function of ${\cal Z}^+$
and seek to understand it as a polynomial in $\lambda$.

This is valid. Following \cite{BNG} we can consider a vector space
$V_{\infty}$ with countable basis $\{ v_0,v_1,v_2,\ldots \}$ on which we
extend the action of $\ls$ by the equations (\ref{rep}).
 For some restriction $\lambda =\Lambda\in {\cal Z}^+$ the
irreducible representation $V_{\Lambda}$ appears as the $\ls$-submodule 
span$\{v_0,v_1,\ldots,v_\Lambda\}$. We know that $W_{\ls}(D)\subset {\cal
Z}(U(\ls))$ so that $\rho_{\Lambda}(W_{\ls}(D))
\propto 1_{\Lambda}$ on the submodule
$V_{\Lambda}$. We define $k(\lambda)[.]:\Ch_*\ra {\cal C}$ by
$\rho_{\lambda}(W_{\ls}(D))= k(\lambda)[D]1_\lambda.$ From
 (\ref{rep}) we see that $k(\lambda)[D]$ must be a finite
polynomial in $\lambda$.
Write $k_{nm}$ for the coefficient of $\lambda^n$ in
$k(\lambda)[.]$ restricted to $\Ch_m$.

The trace over $V_{\Lambda}$ gives  
$W_{\ls,\Lambda}(D)=(\Lambda+1)k(\Lambda)[D]$.
We turn our attention to the weight system that is the coefficient of the
highest order term in $\lambda$ in $\hat{W}_{\ls,\lambda,n}$.
 Thus the weight system we are
investigating is the highest non-vanishing $k_{in}[.]$ (it will turn out
to be $k_{nn}$).
 Our analysis of this weight system is prototypical for all the cases
involved in this paper. Two steps are involved.
\begin{enumerate}
\item{It must be demonstrated that the weight system vanishes on
$I_n^{<n}$. Theorem \ref{cnalg} then implies that such a weight system is in 
the $c$-algebra.}
\item{The exact ${\cal C}$-polynomial in the $\{c_n\}$ is extracted from 
 the values the weight system takes on the basis of $I^n_n$.}
\end{enumerate}

{\bf \underline{Step 1.}}

Take an $n-$CCD D whose graph has $p$ univalent vertices. $W_{\ls}(D)\subset
U(\ls)$ is a sum of products of at most $p$ generators. Equations (\ref{rep})
indicate that $k(\lambda)[D]$ is of maximum order $p$. Recalling that $I_n$ has a
basis of sums of CCDs with no more than $n$ univalent vertices we see that
the highest power of $\lambda$ in $k(\lambda)[D],\ D\in I_n,$ is $n$. Moreover, the maximum 
order on $I^{<n}_n$ is $<n$, so $k_{nn}[I^{<n}_n]=0$. 

\

{\bf \underline{Step 2.}} 

Fix ${\bf P}=\{n_1,\ldots,n_j\}$, an even partition of $n$ and fix $D_{\bf P}$, the planar
representative of $E(\tau_{\bf P})$ (as $k_{nn}$ vanishes on $I^{<n}_n$ it is
constant on the sets $E$). Consider the state sum $W_{\ls,\lambda,n}(D_{\bf
P})$. Only those colourings which colour the legs of $D_{\bf P}$ with the
generator $H$ have non-vanishing weight. This is because the representation
(\ref{rep}) indicates that only the action of $H$ and $X$ can introduce a factor
of $\lambda$ into the matrix elements; moreover, in order that
$\rho_{\lambda}(W_{\ls}) \propto 1_{\lambda}$ any colouring which labelled
some  leg with $X$ would require some other leg to be labelled with a $Y$, the
action of which is independent of $\lambda$. 

Note that $\rho_{\lambda}(H)=\lambda Id\ +\ K$, where $K\in Gl(V_{\lambda})$
is independent of $\lambda$. 
Thus to calculate $k_{nn}[D_{\bf P}]$ we perform the state sum over those
colourings which locate $H$ on all the legs, with the usual weight factors
from
internal vertices, and a factor of $1$ from each univalent vertex. 
This understanding implies:

\bq
k_{nn}[ D_{\{n_1,\ldots,n_j\} }]\ =\ k_{n_1n_1}[ D_{\{n_1\} }]k_{n_2n_2}[ D_{\{n_2\} }]\ldots k_{n_jn_j}[ D_{\{n_j\} }].
\eq

So what is $k_{mm}[ D_{\{m\}} ]$? Recall what $D_{\{m\}}$ looks like. For
example: 
\bq
D_{\{4\}}\  =\  
\Picture{
\put(-1,-1){\line(0,1){0.5}}
\put(-0.4,-1){\line(0,1){0.5}}
\put(0.2,-1){\line(0,1){0.5}}
\put(0.8,-1){\line(0,1){0.5}}
\qbezier(-1,-0.5)(-1.5,1)(-0.1,1)
\qbezier(-0.1,1)(1.3,1)(0.8,-0.5)
\qbezier(-1,-0.5)(-0.7,1)(-0.4,-0.5)
\qbezier(-0.4,-0.5)(-0.1,1)(0.2,-0.5)
\qbezier(0.2,-0.5)(0.5,1)(0.8,-0.5)
\put(-1.2,-1){\vector(1,0){2.3}}
}\ \ .
\eq

\

A colouring of $D_{\{m\}}$ that has a $H$ generator outgoing from some
trivalent
vertex requires that the bracket of the two incoming edges be proportional
to H, if the weight of the colouring is to be non-zero. Recalling that maxima
are labelled with ${\bf b}$, there are two colourings that are non-zero:

\bq
\setlength{\unitlength}{30pt}
\Picture{
\put(-1,-1){\line(0,1){0.5}}
\put(-0.4,-1){\line(0,1){0.5}}
\put(0.2,-1){\line(0,1){0.5}}
\put(0.8,-1){\line(0,1){0.5}}
\qbezier(-1,-0.5)(-1.5,1)(-0.1,1)
\qbezier(-0.1,1)(1.3,1)(0.8,-0.5)
\qbezier(-1,-0.5)(-0.7,1.2)(-0.4,-0.5)
\qbezier(-0.4,-0.5)(-0.1,1.2)(0.2,-0.5)
\qbezier(0.2,-0.5)(0.5,1.2)(0.8,-0.5)
\put(-1.2,-1){\vector(1,0){2.3}}
\put(-1.3,0.75){\mbox{x}}
\put(1,0.75){\mbox{y}}
\put(-1,-0.5){\circle*{0.1}}
\put(0.8,-0.5){\circle*{0.1}}
\put(-0.4,-0.5){\circle*{0.1}}
\put(0.2,-0.5){\circle*{0.1}}
\put(-0.1,1){\circle*{0.1}}
\put(-0.7,0.33){\circle*{0.1}}
\put(-0.1,0.33){\circle*{0.1}}
\put(0.5,0.33){\circle*{0.1}}
\put(-1.02,0.12){\mbox{y}}
\put(-0.55,0.12){\mbox{x}}
\put(0.2,0.12){\mbox{y}}
\put(0.65,0.12){\mbox{x}}
\put(-0.3,-0.35){\mbox{y}}
\put(0,-0.35){\mbox{x}}
\put(-1.3,-0.9){\mbox{h}}
\put(-0.7,-0.9){\mbox{h}}
\put(-0.1,-0.9){\mbox{h}}
\put(0.5,-0.9){\mbox{h}}
}\ \ \ \ ,\ \ \ \
\Picture{
\put(-1,-1){\line(0,1){0.5}}
\put(-0.4,-1){\line(0,1){0.5}}
\put(0.2,-1){\line(0,1){0.5}}
\put(0.8,-1){\line(0,1){0.5}}
\qbezier(-1,-0.5)(-1.5,1)(-0.1,1)
\qbezier(-0.1,1)(1.3,1)(0.8,-0.5)
\qbezier(-1,-0.5)(-0.7,1.2)(-0.4,-0.5)
\qbezier(-0.4,-0.5)(-0.1,1.2)(0.2,-0.5)
\qbezier(0.2,-0.5)(0.5,1.2)(0.8,-0.5)
\put(-1.2,-1){\vector(1,0){2.3}}
\put(-1.3,0.75){\mbox{y}}
\put(1,0.75){\mbox{x}}
\put(-1,-0.5){\circle*{0.1}}
\put(0.8,-0.5){\circle*{0.1}}
\put(-0.4,-0.5){\circle*{0.1}}
\put(0.2,-0.5){\circle*{0.1}}
\put(-0.1,1){\circle*{0.1}}
\put(-0.7,0.33){\circle*{0.1}}
\put(-0.1,0.33){\circle*{0.1}}
\put(0.5,0.33){\circle*{0.1}}
\put(-1.02,0.12){\mbox{x}}
\put(-0.55,0.12){\mbox{y}}
\put(0.2,0.12){\mbox{x}}
\put(0.65,0.12){\mbox{y}}
\put(-0.3,-0.35){\mbox{x}}
\put(0,-0.35){\mbox{y}}
\put(-1.3,-0.9){\mbox{h}}
\put(-0.7,-0.9){\mbox{h}}
\put(-0.1,-0.9){\mbox{h}}
\put(0.5,-0.9){\mbox{h}}}\ \ \ .\
\setlength{\unitlength}{20pt}
\eq

\

\

Both of these colourings have weight $1$: trivalent vertices
contribute factors of $(1)$ in the first colouring, and factors of $(-1)$
in the second, recalling that $m$ is always even. The
maxima all contribute $(1)$. Thus $k_{mm}[D_{\{m\}}]=2$. This implies that
$k_{nn}[D_{\bf P}]=(2)^{\#{\bf P}}$. Glancing at definition (\ref{inv}) indicates
that we have:

\begin{theorem}
\bq
k_{nn}\ =\ \oc_n.
\eq
\rtb
\end{theorem}

\begin{corollary}[MMR]
\bq
JJ(K)(\hb)\ =\ (\sum_n k_{nn} \hb^n)\circ Z_K(K)\ =\ \frac{1}{\tilde{C}(\hb)}.
\eq
\rtb
\end{corollary}

\subsection{$gl(1|1)$}

$\gl$ is a Lie superalgebra on 4 generators: $V_{\overline{0}}=span\{H,G\}, \\
V_{\overline{1}} =span\{Q_+,Q_-\}$. The non-vanishing brackets are the following:
\bq
[G.Q_{\pm}]\ =\ \pm Q_{\pm}\ ,\ \ [Q_+.Q_-]\ =\ H\ .
\eq
There exists an invariant, supersymmetric, non-degenerate form on $\gl$. The
non-vanishing products are:
\bq
<H,G>=-1,\ \ <G,G>=-1,\ \ <Q_+,Q_->=-1.
\eq
This leads to an inverse tensor:
\bq
{\bf b}\ =\ H\otimes H\ -\ G\otimes H\ -\ H\otimes G\ -\ Q_+\otimes Q_-\ +\
Q_-\otimes Q_+.
\eq
$Z(U(\gl))$ is isomorphic to a polynomial ring on two commuting generators
$Z(U(\gl))={\cal C}[c,h]/(ch-hc)$. The isomorphism is $h\mapsto H,\ c\mapsto
G\otimes H\ +\ H\otimes G\ +\ Q_-\otimes Q_+\ -\ Q_+\otimes Q_-$. 

The weight system we are interested in is the deframed 
``universal'' weight system,
$\hat{W}_{\gl,n}:\Ch_n\ra U(\gl)$. Thus we investigate $W_{\gl,n}$
 on the space $I_n$. Consider the state sum for $\gl$.

\begin{lemma}
\label{vanish}
Colourings which label an edge terminating in two trivalent vertices with a
$G$ or an $H$, have weight zero.

{\bf \underline{Proof.}}

First consider turning points. If one edge incident to a maxima
is coloured $G (H)$, then if the other incident edge is not coloured $H (G\,or\,H)$
then the weight of that colouring is zero (consider
{\bf b}). On the other hand, if one edge incident to a minima is coloured
$H(G)$, then if the other incident edge is not coloured $G(G\, or \, H)$ then
the weight of that colouring is zero (consider $\kappa$).

Two points to observe. $H$ annihilates $\gl$ (\ $[H,.]=0$). So if an edge
labelled $H$ is incoming to a trivalent vertex then the weight of that
colouring is zero. Second, observe that $G$ is not in the image of the
bracket. So if an edge labelled $G$ is outgoing from some trivalent vertex,
then the weight of that colouring is zero.

Observe finally that if some internal edge is labelled $G$ or $H$ then taking
care to adjust $G\leftrightarrow H$ at turning points, one of
the above two possibilities is inevitable.
\rtb
\end{lemma}

\begin{lemma}
\bq
W_{\gl}(I^{<n}_{n})\ =\ 0.
\eq
{\bf \underline{Proof.}}

We show that $W_{\gl}$ vanishes on $n-$CCDs whose graphs have $m<n$ univalent
vertices. Choose some such $n-$CCD $D$.
 Call the set of trivalent vertices of $D$ that are linked by an edge
to a univalent vertex $T_1$: there are at most $m$ of these, and at least
$(n-m)$ other vertices. There will be at least one trivalent vertex not in
$T_1$ but linked by an edge to some vertex in $T_1$. This implies that
the diagram for $D$ can be arranged so as to have the following as a
sub-diagram: 

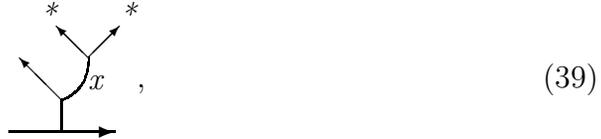
\bq
\Picture{
\thicklines
\put(-1,-1){\vector(1,0){2}}
\thinlines
\put(0,-1){\line(0,1){0.6}}
\put(0,-0.4){\vector(-1,1){0.8}}
\qbezier(0,-0.4)(0.5,-0.2)(0.5,0.4)
\put(0.5,-0.2){\mbox{x}}
\put(0.5,0.4){\vector(1,1){0.6}}
\put(1.1,1){\mbox{*}}
\put(0.5,0.4){\vector(-1,1){0.6}}
\put(-0.4,1){\mbox{*}}
}\ \ ,\
\eq

\

where in the above, the edges marked $(*)$ go to other trivalent
vertices (perhaps by way of some turning points).

We aim to show that every colouring of a diagram which has the above as a
subdiagram has zero weight. Consider the edge marked $x$ in the above. For
the weight to be non-vanishing $x$ must be coloured with a generator in the
image of the bracket. There are two possibilities: $H$ (\ $[Q_+.Q_-]=H$)
and $Q_{\pm}$ (\ $[G.Q_{\pm}]=\pm Q_{\pm}$). If $x$ is coloured $H$ then by the
previous theorem, the weight is zero. Alternatively, if $x$ is coloured by
$Q_{\pm}$ then one of the incoming edges at the top vertex
must have been coloured by $G$. By the previous theorem, the weight of this
colouring also vanishes.
\rtb
\end{lemma}

By the isomorphism mentioned earlier, this weight system may be viewed
$\hat{W}_{\gl}:\Ch_* \ra {\cal C}[c,h]/(ch-hc).$ Then by theorem \ref{cnalg}
 the above  lemma
tells us that the coefficients of particular monomials $c^nh^n$ in
$\hat{W}_{\gl}$ are in the algebra of Alexander-Conway weight systems $\{
c_n\}$. We wish to deduce the exact expressions.

First, recall a property of ``universal'' weight systems $W_V:\Ch_*\ra U(V)$.
\bq
W_V(X.Y)\ =\ W_V(X).W_V(Y),\ \ \mbox{for}\ X,Y\in \Ch_*.
\eq

Choose {\bf P} some partition of $n$. As $W_{\gl}(I^{<n}_n)=0$, $W_{\gl}$ is
constant on the sets $E(\tau_{\bf P})$ so as usual we select the planar
representative $D_{\bf P}$. This satisfies:
\bq
D_{\{n_1,\ldots,n_j\} }\ =\ D_{\{n_1\}}.D_{\{n_2\}}.\ldots.D_{\{n_j\}}.
\eq
Thus,
\bq
W_{\gl}( D_{\{n_1,\ldots,n_j\}})\ =\ W_{\gl}(D_{\{n_1\}})\ldots
W_{\gl}(D_{\{n_j\}}). 
\eq
So, what is $W_{\gl}(D_{\{n\}})$?

\begin{lemma}
\bq
W_{\gl}( D_{\{n\}})\ =\ (-2).H^n
\eq

{\bf \underline{Proof.}}
From lemma \ref{vanish} we know that colourings of $D_{\{n\}}$ with non-vanishing
weights have only $Q_+$ and $Q_-$ colouring internal edges. Labelling maxima
with {\bf b} we observe that there are only two colourings of $D_{\{n\}}$ with
non-vanishing weights. Below we label the edges coloured with $Q_{\pm}$ by a
$\pm$:

\

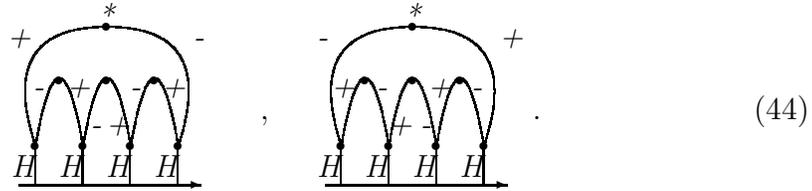
\bq
\setlength{\unitlength}{30pt}
\Picture{
\put(-0.2,1){\mbox{*}}
\put(-1,-1){\line(0,1){0.5}}
\put(-0.4,-1){\line(0,1){0.5}}
\put(0.2,-1){\line(0,1){0.5}}
\put(0.8,-1){\line(0,1){0.5}}
\qbezier(-1,-0.5)(-1.5,1)(-0.1,1)
\qbezier(-0.1,1)(1.3,1)(0.8,-0.5)
\qbezier(-1,-0.5)(-0.7,1.2)(-0.4,-0.5)
\qbezier(-0.4,-0.5)(-0.1,1.2)(0.2,-0.5)
\qbezier(0.2,-0.5)(0.5,1.2)(0.8,-0.5)
\put(-1.2,-1){\vector(1,0){2.3}}
\put(-1.35,0.75){\mbox{+}}
\put(1,0.75){\mbox{-}}
\put(-1,-0.5){\circle*{0.1}}
\put(0.8,-0.5){\circle*{0.1}}
\put(-0.4,-0.5){\circle*{0.1}}
\put(0.2,-0.5){\circle*{0.1}}
\put(-0.1,1){\circle*{0.1}}
\put(-0.7,0.33){\circle*{0.1}}
\put(-0.1,0.33){\circle*{0.1}}
\put(0.5,0.33){\circle*{0.1}}
\put(-1.02,0.12){\mbox{-}}
\put(-0.6,0.12){\mbox{+}}
\put(0.2,0.12){\mbox{-}}
\put(0.6,0.12){\mbox{+}}
\put(-0.3,-0.35){\mbox{-}}
\put(-0.1,-0.35){\mbox{+}}
\put(-1.3,-0.9){\mbox{H}}
\put(-0.7,-0.9){\mbox{H}}
\put(-0.1,-0.9){\mbox{H}}
\put(0.5,-0.9){\mbox{H}}
}\ \ \ \ \ \ ,\ \ \ \ \ \
\Picture{
\put(-0.2,1){\mbox{*}}
\put(-1,-1){\line(0,1){0.5}}
\put(-0.4,-1){\line(0,1){0.5}}
\put(0.2,-1){\line(0,1){0.5}}
\put(0.8,-1){\line(0,1){0.5}}
\qbezier(-1,-0.5)(-1.5,1)(-0.1,1)
\qbezier(-0.1,1)(1.3,1)(0.8,-0.5)
\qbezier(-1,-0.5)(-0.7,1.2)(-0.4,-0.5)
\qbezier(-0.4,-0.5)(-0.1,1.2)(0.2,-0.5)
\qbezier(0.2,-0.5)(0.5,1.2)(0.8,-0.5)
\put(-1.2,-1){\vector(1,0){2.3}}
\put(-1.3,0.75){\mbox{-}}
\put(1,0.75){\mbox{+}}
\put(-1,-0.5){\circle*{0.1}}
\put(0.8,-0.5){\circle*{0.1}}
\put(-0.4,-0.5){\circle*{0.1}}
\put(0.2,-0.5){\circle*{0.1}}
\put(-0.1,1){\circle*{0.1}}
\put(-0.7,0.33){\circle*{0.1}}
\put(-0.1,0.33){\circle*{0.1}}
\put(0.5,0.33){\circle*{0.1}}
\put(-1.12,0.12){\mbox{+}}
\put(-0.55,0.12){\mbox{-}}
\put(0.1,0.12){\mbox{+}}
\put(0.65,0.12){\mbox{-}}
\put(-0.4,-0.35){\mbox{+}}
\put(0,-0.35){\mbox{-}}
\put(-1.3,-0.9){\mbox{H}}
\put(-0.7,-0.9){\mbox{H}}
\put(-0.1,-0.9){\mbox{H}}
\put(0.5,-0.9){\mbox{H}}}\ \ \ .\
\setlength{\unitlength}{20pt}
\eq

\

\

Let us calculate the weights of these colourings. The first has factors -  (-1)
from the maximum labelled $(*)$; (+1) from the other maxima; and (+1) from
each trivalent vertex; the total product is (-1).
 The second has factors - (+1) from
the maximum labelled $(*)$; (-1) from the other maxima; and (-1) from each
trivalent vertex; the total product is again (-1).
\rtb
\end{lemma}

Note that the above differs from the $\ls$ case essentially because in $\gl$ 
 the tensor
{\bf b} changes sign when the $Q_{\pm}$ swap
their order {\bf b}$=\ \ldots + Q_- \otimes Q_+\ -\ Q_+\otimes Q_- \ldots$.
Collecting our results, and comparing with Theorem \ref{cnvalues} we have shown the
following:

\begin{theorem}
\bq
\hat{W}_{\gl,n}(.)\ =\ c_n(.) . H^n\ \ (\subset Z[U(\gl)]).
\eq
\rtb
\end{theorem}

\subsection{ Lie superalgebras of classical type.}
\label{super}

In this section we investigate a class of well-behaved Lie
superalgebras. Recall that an algebra is simple if it has no proper ideals.
\begin{defn}[\cite{Kac}]
A Lie superalgebra $L=L_{\overline{0}}\oplus L_{\overline{1}}$ is of {\it
classical type} if it is simple and if the adjoint representation of
$L_{\zv}\ $ in $\ L_{\ov}$ is completely reducible.
\rtb
\end{defn}

The class of Lie superalgebras of classical type is quite large. The 
classification has been performed in \cite{Kac} and bears many similarities
to the classification of semi-simple Lie algebras, containing several infinite
sequences and several exceptional algebras at distinguished dimensions. More
familiar representatives are $sl(m|n)$ and $osp(m|n)$. $gl(1|1)$ is not of
classical type.

\

%

Hereafter we constrain our analysis to those classical type
Lie superalgebras whose
bodies $L_{\zv}$ are {\it semi-simple Lie algebras}. Some examples of such algebras are the following:

\bq
\begin{array}{|l|l|l|} \hline
L & L_{\zv} & L_{\zv}:L_{\ov} \\ \hline
A(m|m) & A_{m}\oplus A_{m} & sl(m+1)\otimes sl(m+1) \\
B(m|n) & B_m\oplus C_n & so(2m+1)\otimes sp(2n) \\
D(m|n) & D_m \oplus C_n& so(2m)\otimes sp(2n) \\
D(2|1;\alpha) & A_1\oplus A_1\oplus A_1 & sl(2)\otimes sl(2)\otimes sl(2) \\
F(4) & B_3\oplus A_1 & spin(7)\otimes sl(2) \\
G(3) & G_2\oplus A_1 & G_2\otimes sl(2) \\ \hline
\end{array}
\label{list}
\eq
Under $L_{\zv}:L_{\ov}$ we list the representation of $L_{\zv}$ in $L_{\ov}$.
Arkady Vaintrob has pointed out that the following
arguments carry through, largely intact, for 
the algebra $A(m|n)$ when $m\neq n$. In that case, the
$L_{\zv}$ sector has a one-dimensional center.

Take a Cartan subalgebra $H$ of the body of some finite dimensional
Lie superalgebra of classical type $L=L_{\zv}\oplus L_{\ov}$ (\ $H\subset
L_{\zv})$. With this choice, there is a canonical root decomposition 
$L= \oplus_{\alpha\in H^*}L_{\alpha}$ (i.e. if $v\in L_\alpha\ \mbox{then}\ 
[h,v]=\alpha(h)v,\ h\in H$). The set
$\Delta=\{ \alpha\in H^*\ s.t.\ L_{\alpha}\neq \emptyset\}$ is called the
root system of $L$. $\Delta$ partitions $\Delta=\Delta_0\cup \Delta_1$ where
$\Delta_0$ is the root system of $L_{\zv}$ as a Lie algebra and $\Delta_1$ is
the weight system of the representation of $L_{\zv}$ in $L_{\ov}$.

\

\begin{prop}
\label{struct}
Let $L=L_{\zv}\oplus L_{\ov}$ be one of the Lie superalgebras in the list
(\ref{list}) (i.e. at least of classical type with semi-simple body) or a
semi-simple Lie algebra. Write the root decomposition with respect to some
Cartan subalgebra $H$ as $L=\oplus_{\alpha\in H*} L_{\alpha}$. Then:
\begin{enumerate}
\item{ There exists an (unique up to scalar multiple) invariant,
supersymmetric, non-degenerate form on $L$. }
\item{ $L_0 = H$. $H$ is abelian.}
\item{ $dim( L_{\alpha} )=1\ \mbox{when}\ \alpha\neq 0$. }
\item{ $[L_\alpha.L_\beta]\neq 0\ \mbox{iff}\ \alpha,\beta,\alpha+\beta \in
\Delta,$\ in which case\ $[L_\alpha.L_\beta]\subset L_{\alpha+\beta}.$}
\item{ $L_{\alpha}\perp L_{\beta}\ \mbox{if}\ \alpha\neq -\beta.\
<,>|_{L_\alpha \otimes L_{-\alpha}} \neq 0. $}
\item{ $[e_\alpha.e_{-\alpha}]=<e_\alpha,e_{-\alpha}> h_\alpha$, where
$e_\alpha\in L_\alpha,\ e_{-\alpha}\in L_{-\alpha},\ \mbox{and}\ h_\alpha\in
H$\ is defined by $\alpha(h)=<h_\alpha,h>,\ h\in H$. }
\item{ If $\alpha\in \Delta\ (\mbox{resp.} \Delta_0,\Delta_1),\ \mbox{then}\
-\alpha\in \Delta\ (\mbox{resp.} \Delta_0,\Delta_1).$}
\end{enumerate}

{\bf \underline{Proof.}}

This theorem is essentially a restatement of Proposition 5.3 in \cite{Kac}.
The result presented there is much more general, pertaining to the whole
class of classical type Lie superalgebras. Here we have imposed that the
body of L be semi-simple, and we have excluded a number of cases that 
appear as exceptions to particular properties. This lets us strengthen some
of the theorems in an obvious fashion.
\rtb
\end{prop}

\begin{itemize}
\item{ The existence of an appropriate form on $L$ facilitates
the construction of weight systems $W_L$ as per section
\ref{weights}. }
\item{ Choose some linearly independent basis for $\Delta$ over the 
reals $\{ \alpha_1\ldots \alpha_r\}$. Say that $\sum_{i=1}^r c_i \alpha_i \in
H^*$ is positive when $c_k>0$ and $c_j=0,\, 1\leq j<k \leq r$. This partitions
$\Delta=\Delta^+\cup \Delta^-$ into sets of positive and negative roots. We
say that $\alpha>\beta$ when $\alpha-\beta>0$. }
\item{ We can choose canonical generators for the root spaces $L_{\alpha},\
\alpha\neq 0$. Take $\alpha\in \Delta^+$. Item 6 of Theorem \ref{struct} lets
us choose $x_\alpha\in L_\alpha$ and $y_\alpha \in L_{-\alpha}$ such that
$[x_\alpha,y_\alpha]=h_\alpha$, where $h_\alpha\in H$ is defined by the
condition $\alpha(h)=<h_\alpha,h>,\ h\in H$. In this case,
$<x_\alpha,y_\alpha>=1$.}
\end{itemize}

{\bf \underline{Representation theory.}}

Irreducible representations $\rho: L\ra Gl(V)$ of classical Lie superalgebras $L$ satisfy 
a ``highest weight vector'' property \cite{Kac}.
 That is, there exists $v_0\in V$ and  $\lambda\in H^*$
such that,  using the informal
notation $x(v)\equiv \rho_x(v)$,
\begin{enumerate}
\item{$h(v_0)\ =\ \lambda(h)\, v_0,\ h\in H$,}
\item{${x_\alpha}(v_0)\ =\ 0,\ \alpha\in \Delta^+.$}
\end{enumerate}

As with the simpler $sl(2,{\cal C})$ case, it is useful to have a model for
the representations of $L$ in which all irreducible representations
appear as submodules of some formal ``universal'' representation
for appropriate choices of parameters. This lets us
view the class of representations of some element in $U(L)$ 
as a single matrix-valued polynomial in
$\lambda$. Note that when we say that $f$ is a polynomial in $\lambda\in H^*$, we mean that one can introduce some basis in $H^*,\ \{h^i\}$, so that if $\lambda=\lambda_i h^i$ then $f$ is polynomial in the $\lambda_i$. 

Write ${\cal Z}_+\Delta^+$ for the
semi-group  of formal sums of the $\Delta^+$ with coefficients in ${\cal Z}_+$.
We can associate to $k\in {\cal Z}_+\Delta^+$ an element of $U(L)$. 
If $k=
\sum_{\alpha\in \Delta^+}c_{\alpha}[\alpha]\in {\cal Z}_+\Delta^+$ then we write:
\bq
y_k\ =\ \prod_{\alpha\in \Delta^+} y^{c_{\alpha}}_{\alpha},
\eq
where the product is ordered so that if $y_\alpha$ is to the right of
$y_\beta$, $\beta>\alpha$. Define a mapping $i: {\cal Z}_+\Delta^+\ra 
H^*$ by $i(\sum_{\alpha} c_{\alpha} [\alpha])= \sum_{\alpha} c_{\alpha} \alpha$.

Consider the vector space generated by the  following countable basis:
\bq
V_{\infty}\ =\ \{ v_{\bf k}:\ k\in {\cal Z}_+\Delta^+ \}.
\eq
 The action of
$L$ on $V_{\infty}$ is defined as follows, and is a function of
$\lambda\in H^*$:
\begin{eqnarray}
h(v_0) & = & \lambda(h) v_0, \nonumber \\
x_\alpha(v_0) & = & 0, \nonumber \\
y_k ( v_0) &  = & v_k,\ \ k\in {\cal Z}_+\Delta^+.
\end{eqnarray}
These conditions, with the requirement that $V_{\infty}$ represent $U(L)$,
are sufficient to fully determine the representation. Irreducible
representations with highest weights $\Lambda$ appear as submodules when the
specification $\lambda=\Lambda$ is made.

\begin{lemma}
\label{order}
Take some Lie superalgebra from the list (\ref{list}) or a semi-simple Lie
algebra. Then:
\begin{eqnarray}
(1) & h(v_k)\ =\ & (\lambda-i(k))(h)\,v_k,\ \ h\in H,\\
(2) & y_{\alpha}(v_k)\ =\ & \sum_{finite\ sum}d_l v_l, \\
(3) & x_{\alpha}(v_k)\ =\ & \sum_{finite\ sum}c_{l}(\lambda) v_l, 
\end{eqnarray}
where $c_l$ are linear in $\lambda$, and the
$d_l$ are constants. 

{\bf \underline{Proof.}}
\begin{enumerate}
\item{ We know that $[h,y_\alpha]=-\alpha(h)y_{\alpha} (\in U(L))$.
 We seek to calculate
$h(v_k)=h(y_k(v_0))$. Commuting $h$ through the product of $y$ generators produces
our result.}
\item{ Item 4 of Theorem \ref{struct} implies that $y_{\alpha}.y_{\beta}-(-1)^{[y_\alpha][y_\beta]}y_\beta.y_\alpha \subset L_{-(\alpha+\beta)}$ (if $\alpha+\beta$ is a
root) so is proportional to
$y_{\alpha+\beta}$. Left multiplication of some $y_k$ by some $y_\alpha$ is
evaluated by (anti)commuting $y_\alpha$ to it's proper place to produce some sum
over the elements $y_k$. Because only the $y_\alpha$s appear in this process, this is
independent of $\lambda$. }
\item{ Left multiplication of some $y_k$ by some $x_\alpha$ is evaluated
by (anti)commuting the $x_\alpha$ through the product of $y$ generators until it reaches
$v_0$ which it annihilates. Now $x_\alpha.y_\beta-(-1)^{[x_\alpha][y_\beta]} y_\beta.x_\alpha \subset L_{\alpha-\beta}$.
$\alpha-\beta$ may be a positive or negative root, so this process may
create terms containing other $x$, $y$ or $h$ generators. The action of $h$ is
first order in $\lambda$, whence so is the action of $x_{\lambda}$.}
\end{enumerate} 
\rtb
\label{powers}
\end{lemma}

Take $V_{\Lambda}$ some irreducible representation of $L$, and consider some $n-$CCD D. As $W_L(D)\subset
Z(U(L)),\ \rho_{\Lambda}(W_{L}(D))\propto 1_\Lambda$ on the submodule
$V_\Lambda\subset V_\infty$. Define $k(\Lambda)[.]$ by
$\rho_{\Lambda}(W_L(D))=k(\Lambda)[D] 1_\Lambda$. Thus
$W_{L,\Lambda}=k(\Lambda)[D](dim V_{\Lambda})$. Write $k_{nm}$ for the
coefficient of the order $n$ term, when $k(\lambda)[.]$ is restricted to
$\Ch_m$.

\begin{lemma}

\

\begin{enumerate}
\item{$k_{nm}=0\ $if $n>m$.}
\item{$k_{nn}(I^{<n}_n)\ =\ 0.$}
\end{enumerate}

{\bf \underline{Proof.}}

{\bf \underline{(1).}}
We know that $I_n$ has a
basis over $n-$CCDs with not more than $n$ univalent vertices. Thus for $D\in
I_n$, ${W}_{L,n}(D)\subset {\cal Z}[U(L)]$ is a sum of products of at most $n$
generators of $U(L)$. Lemma \ref{powers} indicates that any given generator
can contribute at most one power of $\lambda$ to the matrix: thus 
$k(\lambda)[I_n]$ is of order $\lambda^n$.

\

{\bf \underline{(2).}}
$I^{<n}_n$ has a basis over $n-$CCDs with $<n$ univalent vertices. Thus
$k(\lambda)[.]$
is of order $<n$ in $\lambda$ on $I^{<n}_n$.
\rtb
\end{lemma}

Given item 2 of this lemma, Lemma \ref{cnalg}
 indicates that $k_{nn}[.]$ is in the
algebra of Alexander-Conway weight systems. We wish to derive the exact
expression. Thus we consider the values $k_{nn}$ takes on $I^n_n$. Choose {\bf
P} some partition of $n$, and consider $D_{\bf P}$ the planar representative
of $E(\tau_{\bf P})$. 

Note that in  the state sum $k_{nn}[D_{\bf P}]$ only those colourings of
$D_{\bf P}$ which colour the legs of $D_{\bf P}$ with elements of the Cartan
subalgebra have non-zero weight. This is firstly because Lemma \ref{order}
 indicates
that only the action of the Cartan subalgebra and the raising operators can
contribute powers of $\lambda$ to the matrix elements, and then only first
order. Note further that if any leg is coloured with a raising operator
$X_{\alpha}$ then the weight of that colouring is zero (the matrix
is proportional to the identity so this would require that there be some other
leg coloured with a lowering operator, the action of which is independent of
$\lambda$). 

Observe that $\rho_{\lambda}(h_{\alpha})\ =\ (\lambda(h_{\alpha}) Id + K)$
where $K$ is independent of $\lambda$.
Thus $k_{nn}[D_{\bf P}]$ is calculated from a sum over the weights of
colourings of the graph of $D_{\bf P}$ which colour the legs with elements
of the Cartan subalgebra such that the univalent vertex at the end of an edge
coloured by $h_{\alpha}$ contributes a factor $<\alpha,\lambda>$ to the
weight of that colouring.

This recipe has the immediate consequence that on our planar representative:
\bq
k_{nn}[ D_{\{ n_1,\ldots,n_j\}} ]\ =\ k_{n_1n_1}[ D_{\{n_1\}}]
k_{n_2n_2}[ D_{\{n_2\}}]\ldots k_{n_jn_j}[ D_{\{n_j\}}].
\eq
Thus we need to know $k_{nn}[ D_{\{n\}} ]$.

First consider {\bf b} the tensor that maxima represent. Root vectors appear
in this tensor according to whether they are from the body or the soul. If
$\alpha\in \Delta_0$ we have:
\bq
{\bf b}\ =\ \ldots\ +\ x_{\alpha}\otimes y_{\alpha}\ +\ y_{\alpha}\otimes
x_{\alpha}\ +\ \ldots\ \ \ .\
\eq
Alternatively, if $\alpha\in \Delta_{1}$ we have:
\bq
{\bf b}\ =\ \ldots\ +\ x_{\alpha}\otimes y_{\alpha}\ -\ y_{\alpha}\otimes
x_{\alpha}\ +\ \ldots\ \ \ .\
\label{superb}
\eq

There are two colourings with non-zero weight for each positive 
root. Consider the root $\alpha\in \Delta^+$. In the diagram below, the
label $(+)$ indicates that edge is coloured by $x_{\alpha}$ and the label
$(-)$ indicates that edge is coloured by $y_{\alpha}$.

\bq
\setlength{\unitlength}{30pt}
\Picture{
\put(-0.2,1){\mbox{*}}
\put(-1,-1){\line(0,1){0.5}}
\put(-0.4,-1){\line(0,1){0.5}}
\put(0.2,-1){\line(0,1){0.5}}
\put(0.8,-1){\line(0,1){0.5}}
\qbezier(-1,-0.5)(-1.5,1)(-0.1,1)
\qbezier(-0.1,1)(1.3,1)(0.8,-0.5)
\qbezier(-1,-0.5)(-0.7,1.2)(-0.4,-0.5)
\qbezier(-0.4,-0.5)(-0.1,1.2)(0.2,-0.5)
\qbezier(0.2,-0.5)(0.5,1.2)(0.8,-0.5)
\put(-1.2,-1){\vector(1,0){2.3}}
\put(-1.35,0.75){\mbox{+}}
\put(1,0.75){\mbox{-}}
\put(-1,-0.5){\circle*{0.1}}
\put(0.8,-0.5){\circle*{0.1}}
\put(-0.4,-0.5){\circle*{0.1}}
\put(0.2,-0.5){\circle*{0.1}}
\put(-0.1,1){\circle*{0.1}}
\put(-0.7,0.33){\circle*{0.1}}
\put(-0.1,0.33){\circle*{0.1}}
\put(0.5,0.33){\circle*{0.1}}
\put(-1.02,0.12){\mbox{-}}
\put(-0.6,0.12){\mbox{+}}
\put(0.2,0.12){\mbox{-}}
\put(0.6,0.12){\mbox{+}}
\put(-0.3,-0.35){\mbox{-}}
\put(-0.1,-0.35){\mbox{+}}
\put(-1.4,-0.9){$\mbox{h}_{\alpha}$}
\put(-0.8,-0.9){$\mbox{h}_{\alpha}$}
\put(-0.2,-0.9){$\mbox{h}_{\alpha}$}
\put(0.4,-0.9){$\mbox{h}_{\alpha}$}
}\ \ \ \ \ \ ,\ \ \ \ \ \
\Picture{
\put(-0.2,1){\mbox{*}}
\put(-1,-1){\line(0,1){0.5}}
\put(-0.4,-1){\line(0,1){0.5}}
\put(0.2,-1){\line(0,1){0.5}}
\put(0.8,-1){\line(0,1){0.5}}
\qbezier(-1,-0.5)(-1.5,1)(-0.1,1)
\qbezier(-0.1,1)(1.3,1)(0.8,-0.5)
\qbezier(-1,-0.5)(-0.7,1.2)(-0.4,-0.5)
\qbezier(-0.4,-0.5)(-0.1,1.2)(0.2,-0.5)
\qbezier(0.2,-0.5)(0.5,1.2)(0.8,-0.5)
\put(-1.2,-1){\vector(1,0){2.3}}
\put(-1.3,0.75){\mbox{-}}
\put(1,0.75){\mbox{+}}
\put(-1,-0.5){\circle*{0.1}}
\put(0.8,-0.5){\circle*{0.1}}
\put(-0.4,-0.5){\circle*{0.1}}
\put(0.2,-0.5){\circle*{0.1}}
\put(-0.1,1){\circle*{0.1}}
\put(-0.7,0.33){\circle*{0.1}}
\put(-0.1,0.33){\circle*{0.1}}
\put(0.5,0.33){\circle*{0.1}}
\put(-1.12,0.12){\mbox{+}}
\put(-0.55,0.12){\mbox{-}}
\put(0.1,0.12){\mbox{+}}
\put(0.65,0.12){\mbox{-}}
\put(-0.4,-0.35){\mbox{+}}
\put(0,-0.35){\mbox{-}}
\put(-1.4,-0.9){$\mbox{h}_{\alpha}$}
\put(-0.8,-0.9){$\mbox{h}_{\alpha}$}
\put(-0.2,-0.9){$\mbox{h}_{\alpha}$}
\put(0.4,-0.9){$\mbox{h}_{\alpha}$}}\ \ .
\setlength{\unitlength}{20pt}
\eq

\

\

We must evaluate the weights of these colourings, for $\alpha\in \Delta_0$ and
for $\alpha\in \Delta_1$. Take the first case. Every maxima contributes a
factor $(+1)$ and every trivalent vertex (bracket) contributes a $(+1)$ in the
first colouring $([x_{\alpha},y_{\alpha}]=h_{\alpha})$ and a $(-1)$ in the
second. Keeping in mind the fact that there is alway an even number of
trivalent vertices, the weight from these two colourings is
$2<\lambda,\alpha>^n$.

Take the case $\alpha\in \Delta^1$. The maxima labelled $(*)$ contributes the
opposite sign to rest (consider eqn. (\ref{superb})), so the maxima in both
colourings contribute a factor $(-1)$ to their respective weights. The
brackets always
contribute $(+1)$ (an anticommutator). Thus the weight from these two
colourings is $-2<\lambda,\alpha>^n$.
         
In summary, we have demonstrated that:
\begin{theorem}
\bq
\label{ans}
k_{nn}[D_{ \{n_1,\ldots,n_j\}}]\ =\ \prod_{i=1}^{j}( \sum_{\alpha \in
\Delta^+}(-1)^{[\alpha]}2<\lambda,\alpha>^{n_i} ).
\eq
\rtb
\end{theorem}
 
Collect (formally) the weight systems $\{k_{nn}\}$ as follows:
$K(\hbar)[.]=\sum_{i=0}^{\infty} k_{ii}[.] \hbar^i$. Recall the notations
$c(\hbar)[.]=\sum_{i=0}^{\infty} c_i[.] \hbar^i$ and
$\oc(\hbar)[.]=\sum_{i=0}^{\infty} \oc_i[.]$.

\begin{corollary}
\bq
K(\hbar)[.]\ =\ \prod_{\alpha\in \Delta^+_0} \oc(<\lambda,\alpha>\hbar)
\prod_{\beta \in \Delta_1^+} c(<\lambda,\beta>\hbar).
\eq
{\bf \underline{Proof.}}

Both sides in the above equation vanish on $Ker \phi_n$ and on $I^{<n}_n$. We
must show they coincide for $I_n^n$.
We evaluate $\prod_{\alpha\in \Delta^+_0} c(<\lambda,\alpha>\hbar)
\prod_{\beta \in \Delta_1^+} \oc(<\lambda,\beta>\hbar)
[D_{\{n_1,\ldots,n_j\}}]$. This is equivalent to summing over choices of a
weight system from the product for each component of $D_{\{n_1,\ldots,n_j\}}$
(this requires the fact that $c_n(a.b)=c_{n_1}(a)c_{n_2}(b)$). Recalling
that $c_n(D_{\{n\}})=-2$ and that $\oc_n(D_{\{n\}})=2$ this understanding
exactly reproduces the result (\ref{ans}).
\rtb
\end{corollary}
\begin{corollary}
\bq
K(\hbar)\circ Z_K\ =\ \frac{\prod_{\alpha\in \Delta^+_1}
\tilde{C}(<\lambda,\alpha>\hbar)}{\prod_{\beta\in \Delta^+_0}
\tilde{C}(<\lambda,\beta>\hbar )}.
\eq
\rtb
\end{corollary}

\
	
\noindent{\underbar{\bf Acknowledgements}}
\noindent The author was supported by an Australian Postgraduate Research
Award. He would like to thank Iain Aitchison and Bill Spence for helpful
discussions and a thorough reading of the draft. Thanks also to Arkady
Vaintrob for useful observations on the first draft.

\end{document}